\documentclass[11pt]{article}

\usepackage[T1]{fontenc}
\usepackage[utf8]{inputenc}
\usepackage[english]{babel}
\usepackage[nottoc,notlot,notlof]{tocbibind}
\usepackage{enumerate}
\usepackage{multirow,booktabs}
\usepackage[table]{xcolor}
\usepackage{fullpage}
\usepackage{lastpage}
\usepackage{indentfirst}
\usepackage{fancyhdr}
\usepackage{mathrsfs}
\usepackage{wrapfig}
\usepackage{setspace}
\usepackage{calc}
\usepackage{multicol}
\usepackage{cancel}
\usepackage[retainorgcmds]{IEEEtrantools}
\usepackage[margin=3cm]{geometry}

\usepackage{empheq}
\usepackage{framed}
\usepackage[most]{tcolorbox}
\usepackage{xcolor}
\colorlet{shadecolor}{orange!15}
\usepackage{footnote}

\usepackage{amsmath,amsfonts,amssymb,amscd,amsthm}

\usepackage{bm}

\usepackage[all,2cell]{xy} \UseAllTwocells \SilentMatrices

\usepackage{hyperref}

\usepackage{subfiles}

\newtheorem{teo}{Theorem}[section]
\newtheorem{lem}[teo]{Lemma} 
\newtheorem{cor}[teo]{Corollary}
\newtheorem{prop}[teo]{Proposition} 

\newtheorem{defn}[teo]{Definition} 
\newtheorem{ex}[teo]{Example}

\newtheorem*{claim*}{Claim}

\usepackage{imakeidx}
\makeindex[intoc]

\usepackage{tikz}
\usetikzlibrary{positioning}

\tikzset{onslide/.code args={<#1>#2}{%
  \only<#1>{\pgfkeysalso{#2}} 
}}
\tikzset{temporal/.code args={<#1>#2#3#4}{%
  \temporal<#1>{\pgfkeysalso{#2}}{\pgfkeysalso{#3}}{\pgfkeysalso{#4}} 
the path
}}

  \tikzset{
    invisible/.style={opacity=0},
    visible on/.style={alt={#1{}{invisible}}},
    alt/.code args={<#1>#2#3}{%
      \alt<#1>{\pgfkeysalso{#2}}{\pgfkeysalso{#3}} 
    },
  }
  
\tikzstyle{highlight}=[red,ultra thick]

\usetikzlibrary{mindmap,arrows,shapes.geometric}

\tikzstyle{arrow} = [thick,->,>=stealth]

\usetikzlibrary{positioning}

\tikzset{onslide/.code args={<#1>#2}{%
  \only<#1>{\pgfkeysalso{#2}} 
}}
\tikzset{temporal/.code args={<#1>#2#3#4}{%
  \temporal<#1>{\pgfkeysalso{#2}}{\pgfkeysalso{#3}}{\pgfkeysalso{#4}} 
the path
}}

\tikzstyle{highlight}=[red,ultra thick]

\usetikzlibrary{shapes.geometric, arrows}
\tikzstyle{tp1} = [rectangle, rounded corners, minimum width=2cm, minimum height=1cm,text 
centered, text width=2cm, draw=black, fill=red!30]
\tikzstyle{tp2} = [rectangle, rounded corners, minimum width=2cm, 
minimum height=1cm, text centered, text width=2cm, draw=black, fill=blue!30]
\tikzstyle{tp3} = [rectangle, rounded corners, minimum width=2cm, minimum height=1cm, text 
centered, text width=2cm, draw=black, fill=green!30]
\tikzstyle{tp4} = [rectangle, rounded corners, minimum width=2cm, minimum height=1cm, text 
centered, text width=2cm, draw=black, fill=green!30]
\tikzstyle{tp5} = [rectangle, rounded corners, minimum width=4cm, minimum height=1cm, text 
centered, text width=5cm, draw=black, fill=blue!30]
\tikzstyle{arrow} = [thick,->,>=stealth]

\begin{document}

\author{
  Hugo Rafael de Oliveira Ribeiro\\
  \texttt{hugor@ime.usp.br}
  \and
  Kaique Matias de Andrade Roberto\\
  \texttt{kaique.roberto@usp.br}
  \and
  Hugo Luiz Mariano\\
  \texttt{hugomar@ime.usp.br} \\
  Instituto de Matematica e Estatistica \\
  Universidade de Sao Paulo, Brazil
}

\title{Functorial relationships between multirings and the various abstract theories of quadratic forms}

\maketitle

\begin{abstract}
We provide, explicitly, equivalences and dual equivalences between categories of abstract quadratic forms theories and 
subcategories of multifields and multirings, that will bring new perspectives and methods to the abstract theories of 
quadratic forms in forthcoming papers.

MSC primary 11Exx, secondary 11E81.

\end{abstract}

\section{Introduction}

There are many of abstract theories of quadratic forms. The first ones (abstract Witt rings, quaternionic structures and 
Cordes schemes \cite{marshall1980abstract}) have appeared in the late 70s, by the hands of M. Marshall and C. M. Cordes, with the 
following central target: analyze the existence (or not) of fields with certain properties relating to quadratic forms. In the 
decade of 80's, appears the Marshall's abstract space of orderings (AOS) \cite{marshall1996spaces}: they are important because 
generalize both theory of orderings on fields and the reduced theory of quadratic forms. But only in the early 90's that arise a 
(finitary) first-order theory that generalizes the reduced and non-reduced theory of quadratic forms simultaneously. This theory 
is the special groups of F. Miraglia and M. Dickmann \cite{dickmann2000special}. At that moment, the focus was to look at 
generalizations for the theory of quadratic forms with invertibles coefficients (fields, von Neumman rings, semi-local rings..., 
in general, rings with a good amount of invertibles). In the mid 90's, Marshall generalizes the abstract ordering spaces to 
rings, and called his new theory by ``abstract real spectra'' (ARS), in a first atempt to develop a theory of quadratic forms 
over (general) coefficients on rings. The ring-theoretic case is much more difficult to deal than the field one, the isometry is 
not well behaved and an algebraic counterpart of the abstract real spectra just appears in years 2000, with the real semigroups 
(RS) of Dickmann and Petrovich.

Following the work of professors F. Miraglia and M. Dickmann, through a fruitful and successful partnership between IME-USP and
IMJ-PRG (Paris 6,7), which began in the 1990s, the three authors of this paper continue to expand the boundaries of abstract theories of
quadratic forms, carrying forward the ideas of Dickmann-Miraglia's works, making the IME-USP a center for the development of such 
theories.

All those abstract theories constitute categories that are equivalent, or dually equivalent to full subcategories of each 
other. Also, each one has a particular motivation and advantage. In particular, some of them are categories of first-order 
theories and the corresponding language homomorphisms, thus allowing the application of model-theoretical notions and methods in 
this subject of algebra.


In the present work, we will show that every such theory (and category) can be represented by a category of multirings, that is, 
a ``ring'' with a multivalued addition, a notion introduced in the 1950s by Krasner's works. The notion of 
multiring was joined to the quadratic forms tools by the hands of M. Marshall in last decade (\cite{marshall2006real}). We 
emphasize that multirings by on hand can be described by a first-order theory, and by another allows natural and useful 
generalization of ``multicommutative'' algebra.

\textbf{Overview of the paper}: 

Section 2 contains the preliminary definitions and results on multirings, including 
multicommutative algebra and order theory, needed in the sequel. Most of them (but not all) are already presented in Marshall's 
paper (\cite{marshall2006real}). Section 3 deals with the dual theories of abstract ordering spaces and (reduced) special groups, 
representing them functorially by certain multifields. In section 4 we repeat this process with the abstract theories that deal 
with general coefficients over rings: we represent functorially the dual theories of abstract real spectra and real semigroups 
are represented by Marshall's real reduced multirings, as indicated by Marshall's paper \cite{marshall2006real} and fully 
described here. We emphasize that abstract real spectra, real semigroups (and real reduced multirings) dealt with 
the reduced theory, and a nonreduced approach is not available. At the end of section 4, we connect the new theory of multirings 
and multifields with the most significant theories of quadratic forms. This is (in some way) a new functorial picture: despite 
the Marshall's and Miraglia's observation about these connections, it is the first time that this is made explicit. We summarize 
these functors in the diagram below:

 \SelectTips{eu}{12}
 $$\xymatrix@!=1.6pc{
 \mathcal{RSG}\ar@{^{(}->}[dd]\ar@<.9ex>[rr]\ar[dr] &&
 \mathcal{AOS}^{op}\ar@{^{(}->}'[d][dd]\ar@<.9ex>[ll]\ar@<.9ex>[rr] &&
 \mathcal{MF}_{red}\ar@{^{(}->}[dd]\ar@<.9ex>[ll]\ar[dl]\\
 & \mathcal{SG}\ar@{^{(}->}[dd]\ar@<.9ex>[rr]
 && \mathcal{SMF}\ar@{^{(}->}[dd]\ar@<.9ex>[ll]
 \\
 \mathcal{RS}\ar@<.9ex>[rr]\ar[dr] &&
 \mathcal{ARS}^{op}\ar@<.9ex>[rr]\ar@<.9ex>[ll] &&
 \mathcal{MR}_{red}\ar@<.9ex>[ll]\ar[dl]\\
 & \mathcal{FRS}\ar@<.9ex>[rr]
 && \mathcal{SMR}\ar@<.9ex>[ll]
 }$$
 
 In Section 5, we finalize the paper with indications about how this functorial encoding allows developing concrete steps towards 
{\em non-reduced theory with general coefficients over rings} (\cite{roberto2019nonreduced}) 
{\em and Witt rings associated with those new theories} 
(\cite{ribeiro2019vonneuman}).

\section{Preliminaries}

This section contains the basic definitions and results included for the convenience of the reader. We will deal with multigroups, 
multirings and the analogies of commutative algebra transposed for these new contexts, such as ideals, quotients, orderings and 
etc. For more details, consult \cite{marshall2006real}.

\subsection{Multigroups and Multirings}
\hspace*{\parindent}

Multigroups are a generalization of groups. We can think that a multigroup is a group with a multivalued operation:
\begin{defn}\label{defn:multigroupI}
 A multigroup is a quadruple $(G,\ast,r,1)$, where $G$ is a non-empty set, $\ast:G\times G\rightarrow\mathcal 
P(G)\setminus\{\emptyset\}$ and $r:G\rightarrow G$
 are functions, and $1$ is an element of $G$ satisfying:
 \begin{enumerate}[i -]
  \item If $ z\in x\ast y$ then $x\in z\ast r(y)$ and $y\in r(x)\ast z$.
  \item $y\in 1\ast x$ if and only if $x=y$.
  \item With the convention $x\ast(y\ast z)=\bigcup\limits_{w\in y\ast z}x\ast w$ and 
  $(x\ast y)\ast z=\bigcup\limits_{t\in x\ast y}t\ast z$,
  $$x\ast(y\ast z)=(x\ast y)\ast z\mbox{ for all }x,y,z\in G.$$
  
	A multigroup is said to be \textbf{commutative} if
  \item $x\ast y=y\ast x$ for all $x,y\in G$.
 \end{enumerate}
 Observe that by (i) and (ii), $1\ast x=x\ast 1=\{x\}$ for all $x\in G$. When $a\ast b=\{x\}$ be a unitary set, we just write 
$a\ast b=x$.
\end{defn}

\begin{defn}
 Let $A$ and $B$ be multigroups. A map $f:A\rightarrow B$ is a morphism if for all $a,b,c\in A$:
 \begin{enumerate}[i -]
  \item $c\in a\ast b\Rightarrow f(c)\in f(a)\ast f(b)$;
  \item $f(r(a))=r(f(a))$;
  \item $f(1)=1$.
 \end{enumerate}
\end{defn}


There is another description of multigroups due to M. Marshall\footnote{This is a first-order theory with
axioms of the form $\forall\exists$.}:

\begin{defn}[\cite{marshall2006real}]\label{defn:multigroupII}
 A multigroup is a quadruple $(G,\Pi,r,\mathfrak i)$ where $G$ is a non-empty set, $\Pi$ is a subset of 
$G\times 
G\times G$, $r:G\rightarrow G$ is a
 function and $\mathfrak{i}$ is an element of $G$ satisfying:
 \begin{enumerate}[I -]
  \item If $(x,y,z)\in\Pi$ then $(z,r(y),x)\in\Pi$ and $(r(x),z,y)\in\Pi$.
  \item $(x,\mathfrak{i},y)\in\Pi$ if and only if $x=y$.
  \item If $\exists p\in G$ such that $(u,v,p)\in\Pi$ and $(p,w,x)\in\Pi$ then $\exists q\in G$ such 
that 
$(v,w,q)\in\Pi$ and $(u,q,x)\in\Pi$.
	
	A multigroup is said to be \textbf{commutative} if
  \item $(x,y,z)\in\Pi$ if and only if $(y,x,z)\in\Pi$.
 \end{enumerate}
\end{defn}

 Let $A$ and $B$ be multigroups. A map $f:A\rightarrow B$ is a morphism if for all $a,b,c\in A$:
 \begin{enumerate}[i -]
  \item $(a,b,c)\in\Pi_A\Rightarrow(f(a),f(b),f(c))\in\Pi_B$;
  \item $f(r(a))=r(f(a))$;
  \item $f(1)=1$.
 \end{enumerate}

In fact, these definitions describes the same objects (and arrows), and that connection is established by the following 
lemma:

\begin{lem}[Lemma 1.3 in \cite{marshall2006real}]\label{lemma:1.2}
 For any multigroup $G$ as in the second version, we have:
 \begin{enumerate}[a -]
  \item $r(\mathfrak{i})=\mathfrak{i}$.
  \item $r(r(x))=x$.
  \item $(x,y,z)\in\Pi$ if and only if $(r(y),r(x),r(z))\in\Pi$.
  \item $(\mathfrak{i},x,y)\in\Pi$ if and only if $x=y$.
  \item If $\exists q\in G$ such that $(v,w,q)\in\Pi$ and $(u,q,x)\in\Pi$ then $\exists p\in G$ such 
that 
$(u,v,p)\in\Pi$ and $(p,w,x)\in\Pi$.
  \item For each $a,b\in G$, there exists $c\in G$ such that $(a,b,c)\in\Pi$.
 \end{enumerate}
\end{lem}
 
Now, let $(G,\ast,r,1)$ a multigroup in the sense \ref{defn:multigroupI}. We can define a multigroup 
$(G,\Pi_{\ast},r,\mathfrak i)$ (in the sense \ref{defn:multigroupII}) taking
$\mathfrak i=1$ and $\Pi_{\ast}=\{(a,b,c):c\in a\ast b\}$. The validate of the axioms I,II, III (and IV) 
for $(G,\Pi_{\ast},r,\mathfrak i)$ are direct consequence of axioms i,ii, iii and (iv) in $(G,\ast,r,1)$. 

Conversely, let $(G,\Pi,r,\mathfrak i)$ a multigroup in the sense \ref{defn:multigroupII}. By \ref{lemma:1.2}(f), the function 
$\ast_{\Pi}:A\times A\rightarrow\mathcal P(A)\setminus\{\emptyset\}$, gives by $\ast_{\Pi}(a,b)=a\ast_{\Pi}b:=\{c\in 
G:(a,b,c)\in\Pi\}$ is well-defined. Hence, Let $(G,\ast_{\Pi},1)$ with $1=\mathfrak{i}$. Then, the validate of
the axioms i,ii (and iv) for $(G,\ast_{\Pi},1)$ are direct consequence of I,II, lemma \ref{lemma:1.2}(a) (and IV) for 
$(G,\Pi,r,\mathfrak i)$. For the axiom iii, let $x\in a\ast_{\Pi}(b\ast_{\Pi}c)$. Then $x\in a\ast_{\Pi}q$ for some $q\in 
b\ast_{\Pi}c$. As $(b,c,q)\in\Pi$ and $(a,q,x)\in\Pi$, by \ref{lemma:1.2}(e), there exists $p\in\Pi$ such that $(a,b,p)\in\Pi$ 
and $(p,c,x)\in\Pi$ and then, $x\in p\ast_{\Pi}c$ with $p\in a\ast_{\Pi}b$ imply that
$x\in(a\ast_{\Pi}b)\ast_{\Pi}c$. Finally, let $y\in(a\ast_{\Pi}b)\ast_{\Pi}c$. So $y\in p\ast_{\Pi}c$ for 
some $p\in a\ast_{\Pi}b$, then $(a,b,p)\in\Pi$ and $(p,c,y)\in\Pi$. By III, there exists $q\in\Pi$ such that $(b,c,q)\in\Pi$ and 
$(a,q,y)\in\Pi$. Hence $y\in a\ast_{\Pi}q$ and $q\in b\ast_{\Pi}c$, imply that $y\in a\ast_{\Pi}(b\ast_{\Pi}c)$. Therefore, 
$(G,\ast_{\Pi},1)$ is a multigroup in the sense \ref{defn:multigroupI}.

In fact, this correspondence establishes an isomorphism of (concrete) categories. 

Now, we will define multirings and study this structure in more details. Given a non-empty set $G$ and a multi-operation $\ast \colon G \to 
P(G)$, 
we will use two conventions: if $Z,W\subseteq G$ and $x\in R$, $Z \ast W = \bigcup \{x \ast y: x \in Z,\,y \in W \}$ and 
$Z \ast x = Z \ast \{ x \} = \bigcup \{ z \ast x: z \in Z\}$.

\begin{defn}\label{defn:multiring}
 A multiring is a sextuple $(R,+,\cdot,-,0,1)$ where $R$ is a non-empty set, $+:R\times 
R\rightarrow\mathcal P(R)\setminus\{\emptyset\}$,
 $\cdot:R\times R\rightarrow R$
 and $-:R\rightarrow R$ are functions, $0$ and $1$ are elements of $R$ satisfying:
 \begin{enumerate}[i -]
  \item $(R,+,-,0)$ is a commutative multigroup;
  \item $(R,\cdot,1)$ is a commutative monoid;
  \item $a.0=0$ for all $a\in R$;
  \item If $c\in a+b$, then $c.d\in a.d+b.d$. Or equivalently, $(a+b).d\subseteq a.d+b.d$.
 \end{enumerate}

Note that if $a \in R$, then $0 = 0.a \in (1+ (-1)).a \subseteq 1.a + (-1).a$, thus $(-1). a = -a$.
 
 $R$ is said to be a multidomain if do not have zero divisors; $R$ will be a multifield if every non-zero element of $R$ has 
multiplicative inverse; $R$ is said to be an hyperring if for $a,b,c \in R$, $a(b+c) = ab + ac$. 
\end{defn}

If $F$ be a multifield, then $F$ is an hyperring. Of course, we already have $(a+b)d\subseteq ad+bd$. For the 
other inclusion, if $d=0$, it is done. If $d\ne0$, we have:
 \begin{align*}
 (ad+bd)d^{-1}\subseteq (ad)d^{-1}+(bd)d^{-1}=a+b.
  \end{align*}
  Hence scaling by $d$ we obtain $ad+bd\subseteq(a+b)d$.

\begin{ex}\label{ex:1.3}
$ $
 \begin{enumerate}[a -]
  \item Suppose that $(G,\cdot,1)$ is a group. Defining $a \ast b = \{a \cdot b\}$ and $r(g)=g^{-1}$, 
we have that $(G,\ast,r,1)$ is a multigroup. In this way, every ring, domain and field is a multiring, 
multidomain and multifield, respectively.
  
  \item $Q_2=\{-1,0,1\}$ is multifield with the usual product (in $\mathbb Z$) and the multivalued sum defined by 
relations
  $$\begin{cases}
     0+x=x+0=x,\,\mbox{for every }x\in Q_2 \\
     1+1=1,\,(-1)+(-1)=-1 \\
     1+(-1)=(-1)+1=\{-1,0,1\}
    \end{cases}
  $$
  
  \item In the set $\mathbb R_+$ of non-negative real numbers, we define $a\bigtriangledown 
b=\{c\in\mathbb R_+:|a-b|\le c\le a+b\}$. We have $\mathbb R_+$ with the usual product and $\bigtriangledown$ multivalued sum is 
a multifield, called triangle multifield \cite{viro2010hyperfields}. We denote this multifield by $\mathcal{T}\mathbb R_+$. 
Observe that $\mathcal{T}\mathbb R_+$ is not ``double distributive'': 
$(2\bigtriangledown1)\cdot(2\bigtriangledown1)=[1,3]\cdot[1,3]=[1,9]$ and 
$2\cdot2\bigtriangledown2\cdot1\bigtriangledown1\cdot2\bigtriangledown1\cdot1=
4\bigtriangledown2\bigtriangledown2\bigtriangledown1=[0,9]$.
  
  \item Let $K=\{0,1\}$ with the usual product and the sum defined by relations $x+0=0+x=x$, $x\in K$ 
and 
$1+1=\{0,1\}$. This is a multifield
  called Krasner's multifield \cite{jun2015algebraic}.
  \end{enumerate}
\end{ex}

 Now, another example that generalizes $Q_2=\{-1,0,1\}$. Since this is a new one, we will provide the entire verification that it is a 
multiring:

\begin{ex}[Kaleidoscope]\label{kaleid}
 Let $n\in\mathbb{N}$ and define $X_n=\{-n,...,0,...,n\} \subseteq \mathbb{Z}$. We define the \textbf{$n$-kaleidoscope multiring} by 
$(X_n,+,\cdot,-, 0,1)$, where $- : X_n \to X_n$ is restriction of the  opposite map in$\mathbb{Z}$,  $+:X_n\times 
X_n\rightarrow\mathcal{P}(X_n)\setminus\{\emptyset\}$ is given by the rules:
 $$a+b=\begin{cases}
    \{a\},\,\mbox{ if }\,b\ne-a\mbox{ and }|b|\le|a| \\
    \{b\},\,\mbox{ if }\,b\ne-a\mbox{ and }|a|\le|b| \\
    \{-a,...,0,...,a\}\mbox{ if }b=-a
   \end{cases},$$
and $\cdot:X_n\times X_n\rightarrow X_n$ is is given by the rules:
 $$a\cdot b=\begin{cases}
    \mbox{sgn}(ab)\max\{|a|,|b|\}\mbox{ if }a,b\ne0 \\
    0\mbox{ if }a=0\mbox{ or }b=0
   \end{cases}.$$
   
   .
   In this sense, $X_0=\{0\}$ and $X_1=\{-1,0,1\}=Q_2$. For $X_2$, we have the following "multioperation" table for the sum:
   \begin{center}
\begin{tabular}{|l|l|l|l|l|l|}
\hline
$+$ & $-2$ & $-1$ & $0$ & $1$ & $2$\\
\hline
$-2$ & $\{-2\}$ & $\{-2\}$ & $\{-2\}$ & $\{-2\}$ & $\{-2,-1,0,1,2\}$\\
\hline
$-1$ & $\{-2\}$ & $\{-1\}$ & $\{-1\}$ & $\{-1,0,1\}$ & $\{2\}$\\
\hline
$0$ & $\{-2\}$ & $\{-1\}$ & $\{0\}$ & $\{1\}$ & $\{2\}$\\
\hline
$1$ & $\{-2\}$ & $\{-1,0,1\}$ & $\{1\}$ & $\{1\}$ & $\{2\}$\\
\hline
$2$ & $\{-2,-1,0,1,2\}$ & $\{2\}$ & $\{2\}$ & $\{2\}$ & $\{2\}$ \\
\hline
\end{tabular}
\end{center}
and the following operation table for the product:
\begin{center}
\begin{tabular}{|l|l|l|l|l|l|}
\hline
$\cdot$ & $-2$ & $-1$ & $0$ & $1$ & $2$\\
\hline
$-2$ & $2$ & $2$ & $0$ & $-2$ & $-2$\\
\hline
$-1$ & $2$ & $1$ & $0$ & $-1$ & $-2$\\
\hline
$0$ & $0$ & $0$ & $0$ & $0$ & $0$\\
\hline
$1$ & $-2$ & $-1$ & $0$ & $1$ & $2$\\
\hline
$2$ & $-2$ & $-2$ & $0$ & $2$ & $2$ \\
\hline
\end{tabular}
\end{center}
   
   Clearly $(X_n, \cdot,1)$ is a commutative monoid and $a\cdot0=0$ for all $a\in X_n$.
   
Now, we will verify that $(X_n,+,\cdot, -, 0,1)$ is a multiring. 
   \begin{enumerate}[i -]
    \item By construction, $a+b=b+a$, $a+0=\{a\}$ and $0\in a-a$ for all $a,b\in X_n$.
    
    \item $d\in a+b\Leftrightarrow b\in d-a$: We divide the proof in cases. Let $a\ne-b$ and suppose without loss of generality that $|a|<|b|$. 
Thus $a+b=\{b\}$. Hence $d\in a+b$ implies $d=b$. So $b\in b-a=\{b\}$. By symmetry, the same proof applies to the implication $b\in 
d-a\Rightarrow d\in a+b$. The case $|b|=|a|$ is immediate.
    
    \item $(a+b)+c=a+(b+c)$: Again we divide in cases. We suppose without loss of generality that $a,b,c\ne0$. 
If $a\ne-b$, $b\ne-c$, and $|a|\le |b|\le |c|$,
$$(a+b)+c=a+(b+c)=\{c\}.$$
Similarly, $(a+b)+c=a+(b+c)$ for the cases $|a|\le |c|\le |b|$, $|b|\le |a|\le |c|$, $|b|\le |c|\le |a|$,
$|c|\le |a|\le |b|$ and $|c|\le |b|\le |a|$ (under the hypothesis $a\ne-b$, $b\ne-c$).
    
    Now let $a=-b$. We want to prove that $(a-a)+c=a+(-a+c)$. If $|a|\le|c|$, 
    $$(a-a)+c=X_a+c=\{c\}\mbox{ and }a+(-a+c)=a+c=\{c\}.$$ 
    If $|c|<|a|$, then 
    $$(a-a)+c=X_a+c=X_a\mbox{ and }a+(-a+c)=a-a=X_a$$ 
    The case $b=-c$ is analogous. 
    
    
    \item $d(a+b)\subseteq da+db$: If $d=0$ there is nothing to prove. Let $d\ne0$. If $a\ne-b$, suppose without loss of 
generality that $|a|<|b|$. Then $a+b=\{b\}$ and $d(a+b)=\{db\}=db+db$.
    
    Now let $a=-b$. We have two cases:
    \begin{enumerate}
     \item $|d|\le|a|$: since $da=\mbox{sgn}(da)|a|$, we have $da-da=X_{da}=X_a$ and $d(a-a)=dX_a\subseteq X_a$.
     
     \item $|d|>|a|$: since $da=\mbox{sgn}(da)|d|$, we have $da-da=X_{da}=X_d$ and $d(a-a)=dX_a\subseteq X_d$.
    \end{enumerate}
   \end{enumerate}
   
   Thus $X_n$ is a multiring.
\end{ex}

 Now, another example that generalizes $K=\{0,1\}$. Since this is a new one, we will provide the entire verification that it is a {\em  
multifield}:

\begin{ex}[H-multifield]\label{H-multi}
Let $p\ge1$ be a prime integer and $H_p:=\{0,1,...,p-1\} \subseteq \mathbb{N}$. 
Now, define the binary multioperation and operation in $H_p$ as follow:
\begin{align*}
 a+b&=
 \begin{cases}H_p\mbox{ if }a=b,\,a,b\ne0 \\ \{a,b\} \mbox{ if }a\ne b,\,a,b\ne0 \\ \{a\} \mbox{ if }b=0 \\ \{b\}\mbox{ if }a=0 \end{cases} \\
 a\cdot b&=k\mbox{ where }0\le k<p\mbox{ and }k\equiv ab\mbox{ mod p}.
\end{align*}
By a similar argument used in example \ref{kaleid} we obtain that $(H_p,+,\cdot,-, 0,1)$ is a multifield such that for all $a\in H_p$, $-a=a$. 
For example, considering $H_3=\{0,1,2\}$, using the above rules we obtain these tables
\begin{center} 
 \begin{tabular}{|l|l|l|l|l|l|}
\hline
$+$ & $0$ & $1$ & $2$\\
\hline
$0$ & $\{0\}$ & $\{1\}$ & $\{2\}$ \\
\hline
$1$ & $\{1\}$ & $\{0,1,2\}$ & $\{1,2\}$\\
\hline
$2$ & $\{2\}$ & $\{1,2\}$ & $\{0,1,2\}$ \\
\hline
\end{tabular}
\end{center}

\begin{center} 
 \begin{tabular}{|l|l|l|l|l|l|}
\hline
$\cdot$ & $0$ & $1$ & $2$\\
\hline
$0$ & $0$ & $0$ & $0$ \\
\hline
$1$ & $0$ & $1$ & $2$\\
\hline
$2$ & $0$ & $2$ & $1$ \\
\hline
\end{tabular}
\end{center}
In fact, these $H_p$ is a kind of generalization of $K$, in the sense that $H_2=K$. Now, we proceed with the verifications. Clearly $(H_p, 
\cdot,1)$ is a commutative monoid and $(H_p\setminus\{0\},\cdot,1)$ is an abelian group. Moreover, $a\cdot0=0$ for all $a\in H_p$.
 
 Now, we will verify that $(H_p,+,\cdot,-, 0,1)$ is a multiring:
 \begin{enumerate}[i -]
    \item By construction, $a+b=b+a$, $a+0=\{a\}$ and $0\in a-a = a+a$ for all $a,b\in H_p$.
    
    \item $d\in a+b\Leftrightarrow b\in d-a = d+a$: We split the proof in cases. The case where $a = 0$ or $b=0$ is immediate, so we treat the 
case $a, b  \neq 0$. If $a =b$, then $a+a =  H_p$ and  $d \in a+a \ {\rm{iff}} \ a \in d +a$.  Let $a\ne b$, thus $a+b=\{a,b\}$, and $d\in 
a+b\Leftrightarrow$ $d=a$ or $d=b$. If $d=a$, then $d=a\in b-b=b+b$, and if $d=b$, we obtain $d=b\in a-a=a+a$, proving the desired property.
    
    
    \item $(a+b)+c=a+(b+c)$: Again we divide in cases. We suppose without loss of generality that $a,b,c\ne0$. 
If $a\ne b$, $b\ne c$ and $a\ne c$, then 
$$(a+b)+c=a+(b+c)=\{a,b,c\}.$$
    
    Now let $a= b$ (the case $b=c$ is similar). We want to prove that $(a+a)+c=a+(a+c)$. 
    $$(a+a)+c=H_p+c=H_p\mbox{ and }a+(a+c)=a+\{a, c\}=H_p.$$ 
    
    
    \item $d(a+b)\subseteq da+db$: If $d=0$ there is nothing to prove. Let $d\ne0$. If $a=0$ or $b=0$, again, there is nothing to prove. Let 
$a,b,d\ne0$. If $a=b$, then
    $$d(a+b)=d(a+a)=dH_p=H_p=da+da = da+db.$$
    If $a\ne b$, then $da \neq db$ and 
    $$d(a+b)=d\{a,b\}=\{da,db\}=da+db.$$
   \end{enumerate}
   Thus $H_p$ is a multifield.
\end{ex}

 Now, we treat about morphisms:

\begin{defn}\label{defn:morphism}
 Let $A$ and $B$ multirings. A map $f:A\rightarrow B$ is a morphism if for all $a,b,c\in A$:
  \begin{enumerate}[i -]
  \item $c\in a+b\Rightarrow f(c)\in f(a)+f(b)$;
  \item $f(-a)=-f(a)$;
  \item $f(0)=0$;
  \item $f(ab)=f(a)f(b)$;
  \item $f(1)=1$.
 \end{enumerate}

\end{defn}

 For multirings, there are types of morphisms that can be considered. Let $f \colon A \to B$ a multiring morphism.
 \begin{itemize}
    \item $f$ is a \textbf{strong morphism} if for all $a,b, c\in A$, if $f(c) \in f(a) + f(b)$, then there exist $a',b',c' \in A $ with $f(a') = 
f(a),f(b') = f(b), f(c') = f(c)$ such that $c' \in a' + b'$.
    \item $f$ is an \textbf{ideal morphism} if for all $a,b, c\in A$, if $f(c) \in f(a) + f(b)$, then exists $c' \in A $ with $f(c') = f(c)$ such 
that $c' \in a + b$. In other words, $f(a + b) = (f(a) + f(b)) \cap \mbox{Im}(f)$
     \item We say that $f$ is a \textbf{full morphism} if it is a strong morphism for all $a,b\in A$ and all $d\in B$,
     $$d\in f(a)+f(b) \Rightarrow \mbox{ exists } c \in a + b \mbox{ such that }d = f(c).$$
     In other words, $f(a+b)=f(a)+f(b)$.
     \item We say that $f$ is a \textbf{strong embedding} if $f$ is injective and it is a strong morphism. In this case, $A$ is a 
\textbf{submultiring} of $B$ if $A\subseteq B$ and the canonical inclusion $\iota:A\hookrightarrow B$ is a strong embedding.
     \item We say that $f$ is a \textbf{full embedding} if it is a strong embedding and a full morphism\footnote{There is no consensus on the 
definition "submultiring": here we do  adopted one of intermediary strength  that coincides with the notion of substructure in relational 
structures; in \cite{marshall2006real}, submultiring means an inclusion of multirings that is strong {\em and} full.}.
 \end{itemize}
 
 The following diagram
illustrates the diferent notions of morphisms and their relations.
 $$\mbox{Full Morphism}\Rightarrow\mbox{Ideal Morphism}\Rightarrow\mbox{Strong Morphism}$$
$$\mbox{Full Embedding}\Rightarrow\mbox{Strong Embedding}\Leftrightarrow
\mbox{Ideal Embedding}$$
 
 If $f \colon A \to B$ is multiring morphism with $B$ ring, then $f$ is a full morphism. Furthermore, an isomorphism in multiring category is a 
full embedding.
 
 The category of multifields (respectively multirings) and their morphisms will be 
denoted by $\mathcal{MF}$ (respectively $\mathcal{MR}$).

Some of the properties of rings morphisms are not extended to multirings morphisms. In the sequel there are some examples and 
counterexamples:

 \begin{ex}\label{ex:2.1}
  $ $
  \begin{enumerate}[a -]
   \item Let $f:\mathbb R\rightarrow Q_2$ be $f(x)=\mbox{sgn}(x)$, the signal function (with convention that $\mbox{sgn}(0)=0$). 
$f$ is a multiring morphism, but $f$ is not injective and $\mbox{Ker}f=\{0\}$. Furthermore, $f$ is a strong morphim that it is not an ideal 
morphism, and therefore, also it is not a full morphism.
   
   \item Consider $\mathbb F_2=\{0,1\}$ (the field with two elements) and let $K=\{0,1\}$ be the multifield with the structure $1+0=0+1=\{1\}$, 
$1+1=\{0,1\}$, and product as usual (in $\mathbb Z$). Consider the inclusion $i:\mathbb F_2\rightarrow K$, i.e, $i(0)=0$ and $i(1)=1$. Then $i$ 
is an embedding that is not strong (and hence, not ideal neither full):
   $$i(1+1)=\{0\}\ne\{0,1\}=i(1)+i(1),$$
   hence $i(1+1)\ne(i(1)+i(1))$.
   
   \item The inclusions functions $Q_2\hookrightarrow\mathbb R$ and $\mathcal T\mathbb R_+\hookrightarrow\mathbb R$ are not multiring morphisms.
   
if $a\ne0$, $a^2=1$. For example, considering $H_2=\{0,1,2\}$, using the above rules we obtain these tables


\item Consider the set of multifields $H_p, p \geq 1$, in  Example \ref{H-multi}, where $H_2 =K$.  Then inclusion morphism $j: H_2\rightarrow 
H_3$ given by the rule $j(0)=0$, $j(1)=1$.Then $j$ is an ideal embedding that is not full:
$$j(1+1)=\{0,1\}\ne H_3=j(1)+j(1).$$

\item  Consider the diagonal morphism $\Delta \colon Q_2 \to Q_2 \times Q_2$, given by the rule $\Delta(x)=(x,x)$. Then $\Delta$ is a strong 
embedding that is not full: 
$$(1,-1) \in \Delta(1) + \Delta(-1), \ {\rm but} \ (1,-1) \notin \mbox{Im}(\Delta).$$
  \end{enumerate}
 \end{ex}
 
 \subsection{Commutative Multialgebra}
 Here, we will extend some terminology of commutative algebra to multirings and multifields. As expected, many concepts such 
that ideals, fractions, and localizations have a natural generalization for multirings. We register that some of the results below 
seem to be new, or at least, unpublished. For other details about this subject, the reader can consult 
\cite{ramaruban2014commutative}.

 \begin{defn}\label{defn:ideal}
  An ideal of a multiring $A$ is a non-empty subset of $A$ such that $\mathfrak{a}+\mathfrak{a}\subseteq\mathfrak{a}$ and 
$A\mathfrak{a}= \mathfrak{a}$. An ideal $\mathfrak{p}$ of $A$ is said to be prime if $1\notin\mathfrak{p}$ and $ab\in\mathfrak{p}\Rightarrow 
a\in\mathfrak{p}$ or $b\in\mathfrak{p}$. An ideal $\mathfrak{m}$ is maximal if for all ideals $\mathfrak{a}$ with 
$\mathfrak{m}\subseteq\mathfrak{a}\subseteq 
A\Rightarrow$, then $\mathfrak{a}=\mathfrak{m}$ or $\mathfrak{a}=A$. We will denote $\mbox{Spec}(A)=\{\mathfrak{p}\subseteq 
A:\mathfrak{p}\mbox{ is a prime ideal}\}$.
 \end{defn}
 
If $\mathfrak{a}$ is an ideal of $A$, note that $0 \in \mathfrak{a}$ and  $-\mathfrak{a} \subseteq \mathfrak{a}$.
 
 With the notion of ideal, we can define some new multirings structures with the language of commutative algebra in mind:
 \begin{defn}\label{terminology}
  $ $
  \begin{enumerate}[a -]
   \item If $\{A_i\}_{i\in I}$ is a family of multirings, then the product $\Pi_{i\in I}A_i$ is a multiring in the natural 
(component wise) way.
   
   \item Let $\mathfrak{a}\subseteq A$ be an ideal. Elements of $A/\mathfrak{a}$ are cosets $\overline a=a+\mathfrak{a}$, $a\in A$. More 
explicitly,
   $$a\equiv b\mbox{ mod }I\mbox{ if and only if }b\in\overline a\mbox{, if and only if }(b-a) \cap \mathfrak a \neq \emptyset.$$
This is the multialgebra analogous of the usual congruece relation in commutative algebra. We define a multiring structure on $A/\mathfrak{a}$ by 
$\overline a+\overline b=\{\overline c:c\in a+b\}$, 
$-\overline a=\overline{-a}$, the zero and the unit element of $A/\mathfrak{a}$ are $0=\overline{0}$ and $1=\overline{1}$ 
respectively and multiplication on $A/\mathfrak{a}$ is defined by $\overline a\overline b=\overline{ab}$. Note that if $\overline{c} \in 
\overline{a} + \overline{b}$,
then exists $c' \in a + b$ such that $\overline{c'} = \overline{c}$. The natural arrow
$\pi: A \rightarrow A /\mathfrak{a}$ is a strong morphism and as in the ring case it is easily proved that given another multiring 
morphism $f: A \rightarrow B$ with $f(\mathfrak{a}) = \{0\}$, there is a unique morphism $\overline{f}: A / \mathfrak{a} 
\rightarrow B$ such that $f = \overline{f} \circ \pi$.
   
   \item Let $S$ be a multiplicative set in $A$. Elements of $S^{-1}A$ have the form $a/s$, $a\in A$, $s\in S$, $a/s=b/t$ if 
and only if  $atu=bsu$ for some $u\in S$. $0=0/1$, $1=1/1$ and the operations are defined by $(a/s)\cdot(b/t)=ab/st$, and $c/v\in 
a/s+b/t$ if and only if $cstv\in atuv+bsuv$ for some $v\in S$. The natural arrow $\rho \colon A\hookrightarrow S^{-1}A$ is a strong morphism and 
given
a multirng morphism $f \colon A \to B$ with $f(S) \subseteq B^{\times}$, then exists a unique morphism $\overline{f} \colon S^{-1} A \to B$ such 
that 
$f = \overline{f} \circ \rho$.
   
   \item If $D$ is a multidomain, we define the multifield of fractions $\mbox{ff}(D):=(D\setminus\{0\})^{-1}D$.
  \end{enumerate}
 \end{defn}

 Let $X$ be a  subset of a multiring $A$. We define the \textbf{ideal generated by} $X$ as $\langle 
X \rangle:=\bigcap\{\mathfrak{a}\subseteq A : X \subseteq\mathfrak{a}, \mathfrak{a} \ \mbox{ is an ideal} \}$. If $X \neq \emptyset$, it can 
easily checked that
 $\langle X \rangle=  \bigcup \{\lambda_1x_1 + \cdots + \lambda_n x_n \colon  n \geq 1, \lambda_i \in A, x_i \in X, \mbox{ for all } i=1, \ldots, 
n\}$. In particular $$\langle a \rangle = \sum Aa :=
 \left\lbrace\sum^n_{j=1}\lambda_ja:\lambda_1,...,\lambda_n\in A,\,n\ge1.\right\rbrace.$$
 
If $A$ is a hyper-ring, i.e. if it  satisfies also the second-half distributive, then $\sum Aa=Aa$.

 \begin{prop}\label{lem:iso}
  Let $A$ and $B$ be multirings and $\varphi: A \rightarrow B$ a surjective morphism. Consider $\overline{\varphi}: A /\mbox{Ker} 
(\varphi) \rightarrow B$ the induced morphism. Then the following are equivalent:
  
  \begin{enumerate}[i)]
   \item $\varphi$ is a strong morphism and if $\varphi(a) = \varphi(a')$ for $a,a' \in A$, then $(a - a') \cap \mbox{ker}(\varphi) \neq 
\emptyset$.	
   
   \item $\varphi$ is an ideal morphism.
   
   \item $\overline{\varphi}$ is an isomorphism.
  \end{enumerate}
 \end{prop}
 \begin{proof}
  $i) \Rightarrow ii)$: Assume that $\varphi(a) \in \varphi(b) + \varphi(c)$. Since $\varphi$ is a strong morphism, exists $a',b',c' \in A$ with
  $\varphi(a') = \varphi(a), \varphi(b') = \varphi(b), \varphi(c') = \varphi(c)$ such that $a' \in b' + c'$. By hypothesis, exists $b' \in b + i$ 
and $c' \in c + j$ 
  such that $i,j \in \mbox{ker}(\varphi)$. Then $a' \in b' + c' \subseteq (b + c) + (i + j)$ and so exists $x \in i + j \subseteq 
\mbox{ker}(\varphi)$ such that 
  $a' \in b + c + x$. Thus exist $a'' \in a' - x$ with $a'' \in b + c$ and note that $\varphi(a'') = \varphi(a') = \varphi(a)$.
 
  $ii) \Rightarrow iii)$: Let $a,b \in A$ such that $\varphi(a) = \overline{\varphi}(\overline{a}) = \overline{\varphi}(\overline{b}) = 
\varphi(b)$. 
  By hypothesis exist $x \in a - b$ such that $x \in \mbox{ker}(\varphi)$ and so $\overline{a} = \overline{b}$ in $A / \mbox{ker}(\varphi)$, 
  proving the injectivity of $\overline{\varphi}$. Since $\varphi$ is a strong morphism, if $\overline{\varphi}(\overline{a}) \in 
\overline{\varphi}(\overline{b}) + 
\overline{\varphi}(\overline{c})$, then exists $a',b',c' \in A$ with $\varphi(a') = \varphi(a), \varphi(b') = \varphi(b), \varphi(c') = 
\varphi(c)$ such that
$a' \in b' + c'$. By hypothesis, it is easy to see that $\overline{a'} = \overline{a}, \overline{b'} = \overline{b}, \overline{c'} = 
\overline{c}$ and so
$\overline{a} \in \overline{b} + \overline{c}$ in $A/\mbox{ker}(\varphi)$. Thus $\overline{\varphi}$ is an isomorphism.
  
  $iii) \Rightarrow i)$: Assume that $\varphi(a) = \varphi(a')$ for $a,a' \in A$. Then $\overline\varphi(\overline a)=\overline\varphi(\overline 
a')$ and hence $\overline a=\overline{a'}$, which means that $(a-a')\cap\mbox{ker}(\varphi)\ne\emptyset$.
  
  Therefore $0 = \varphi(0) \in \varphi(a) - \varphi(a')$, and by hypothesis there exist $i \in a - a'$ such that
  $\varphi(i) = \varphi(0) = 0$. On the other hand, we have $\varphi = \overline{\varphi} \circ \pi$, where $\pi \colon A \to A/ 
\mbox{ker}(\varphi)$. Then $\varphi$ is a composition
  of strong morphisms and so $\varphi$ is strong itself.
 \end{proof}
 
 \begin{teo}[Isomorphism Theorem]\label{teo:iso}
  Let $A$ and $B$ be multirings and $\varphi:A\rightarrow B$ an ideal morphism. Then $\mbox{Im}(\varphi)$ is a multiring (contained in $B$) with 
the structure induced by the domain $A$, and the induced morphism $\overline{\varphi}: A/\mbox{Ker}(\varphi) \rightarrow \mbox{Im}(\varphi)$ is 
an isomorphism.
 \end{teo}
 \begin{proof} 
  By the previous proposition, it is enough to prove that $\mbox{Im} (\varphi)$ is a multiring and this is accomplished by proving the 
associativity property for $\mbox{Im}(\varphi)$. Assume that $\varphi(x) \in \varphi(p) + \varphi(w)$ with $\varphi(p) \in 
\varphi(u) + \varphi(v)$. Since $\varphi$ is an ideal morphism, exists $x' \in p + w$ and $p' \in u + v$ such that $\varphi(x') = 
\varphi(x)$ and $\varphi(p') = \varphi(p)$. Then, by the same argument as the previous lemma, it should exist $i \in \mbox{Ker} 
(\varphi)$ such that $p \in i + p'$. Then $p \in i + (u + v) \subseteq (i + u) + v$ and thus exist $u' \in i + u$ such that $p 
\in u' + v$. Then exist $q \in v + w$ with $x \in u' + q$. Thus $\varphi(x) \in \varphi(u') + \varphi(q) = \varphi(u) + 
\varphi(q)$ and $\varphi(q) \in \varphi(v) + \varphi(w)$.
 \end{proof}
 
\begin{lem}
Let $A$ be a multiring and $\mathfrak{a}$ and ideal. Consider $\pi:A \to A/\mathfrak{a}$ the canonical projection. Then it 
induces a one-to-one correspondence between the ideals of $A$ that contains $\mathfrak{a}$ and the ideals of $A/\mathfrak{a}$.
\end{lem} 
 \begin{proof}
    The proof is the same as the ring case.
 \end{proof}
 
 \begin{lem}\label{lemma:1.1}
  Let $A$ be a multiring. Then:
  \begin{enumerate}[a -]
   \item an ideal $\mathfrak{p}$ of $A$ is prime if and only if $A/\mathfrak{p}$ is a multidomain.
   \item An ideal $\mathfrak{m}$ is maximal if and only if for all $a \neq 0$ in $A/\mathfrak m$, exists $t_1, \ldots, t_n$ such 
that $1 \in at_1 + \cdots + at_n$. In particular, maximal ideals are prime and if $A$ is a hyperring, an ideal $\mathfrak{m}$ is 
maximal if and only if $A/\mathfrak{m}$ is a multifield.
  \end{enumerate}
 \end{lem}
 \begin{proof}
   \begin{enumerate}[a -]
      \item The same of the ring case.
      \item $\Rightarrow$: Let $\overline{a} \in A/\mathfrak{m}$ non-zero, that is, $a \notin \mathfrak{m}$. Since $m$ is 
maximal, 
the ideal $I = \bigcup \{m + at_1 + \cdots + at_n \colon n \geq 1 \mbox{ and } t_i \in A\}$ generated by $\mathfrak{m} \cup 
\{a\}$ is improper. Then exists $m \in \mathfrak{m}$ and $t_1, \cdots, t_n \in A$ such that $1 \in m + at_1 + \cdots + at_n$
and so $\overline{1} \in \overline{a} \overline{t_1} + \cdots + \overline{a} \overline{t_n}$. 

$\Leftarrow$: Let $a \notin 
\mathfrak{m}$. By the property valid in $A/\mathfrak{m}$, exists $m \in \mathfrak{m}$ and $t_1, \cdots, t_n \in A$ such that $1 
\in m + at_1 + \cdots + at_n$ and so the ideal generated by $\mathfrak{m} \cup\{a\}$ is improper. Then $\mathfrak{m}$ is maximal.
   \end{enumerate}
 \end{proof} 

 \begin{prop}
 $ $
  \begin{enumerate}[a -]
   \item Let $A$ be a multiring, $I\subseteq A$ an ideal and $S\subseteq A$ be a multiplicative subset of $A$. Then 
$(S/I)^{-1}A/I \cong S^{-1}A / S^{-1}I.$

   \item Let $\{A_i\}_{i\in I}$ be a family of multirings and $\mathfrak a_i\subseteq A_i$ be an ideal of $A_i$ for every $i\in 
I$. Then $$\prod_{i\in I}A_i/\mathfrak a_i\cong\prod_{i\in I}A_i/\prod_{i\in I}\mathfrak a_i.$$
  \end{enumerate}
 \end{prop}
 \begin{proof}
  For the item (a), consider the morphism $f:S^{-1} A \rightarrow (S/I)^{-1} A/I $ given by  $f(a/s)=\overline a/\overline s$ and 
apply the theorem \ref{teo:iso}. For the item (b), the same strategy holds with the morphism 
$g:\prod_{i\in I}A_i\rightarrow\prod_{i\in I}A_i/\prod_{i\in I}\mathfrak a_i$ given by $g(a_i)_{i\in I}=(\overline a_i)_{i\in I}$.
 \end{proof}

$ $

 Now, we present a construction that will be used several times below:
 
 \begin{defn}\label{defn:strangeloc}
 Fix a multiring $A$ and a multiplicative subset $S$ of $A$ such that $1\in S$. Define an equivalence relation $\sim$ 
on $A$ by $a\sim b$ if and only if $as=bt$ for some $s,t\in S$. Denote by $\overline a$ the equivalence class of 
$a$ and set $A/_mS=\{\overline a:a\in A\}$. Then, we define in agreement with Marshall's notation, $\overline a+\overline 
b=\{\overline c:cv\in as+bt,\,\mbox{for some }s,t,v\in S\}$, $-\overline a=\overline{-a}$, and 
$\overline{a}\overline{b}=\overline{ab}$.\footnote{
If $0\in S$, then $A/_mS=\{\overline0\}$. We keep this trivial case available  for precaution in definition \ref{defn:strangeloc}. One of the 
possible aspects to be investigated after the considerations of this paper is the creation of a correspondence $(F,T)\mapsto F/_m\dot T$, for 
suitable category which objects are pairs $(F,T)$, where $F$ is a field and $T\subseteq F$. It is desirable for a such category of pairs to 
provide good categorical properties (for example, closeness by finite products and colimits), and for a such correspondence, to provide good 
functorial properties. In this hypothetical scenario, the trivial Marshall's quotient may play the role of zero object.
}
 \end{defn}
 
 In particular, the canonical projection $\pi:A\rightarrow A/_mS$ is a strong morphism.
 \begin{prop}
  Let $A,B$ be a multiring and $S\subseteq A$ a multiplicative subset of $A$. Then for every morphism 
$f:A\rightarrow B$ such that $f[S]=\{1\}$, there exist a unique morphism $\tilde f:A/_mS\rightarrow B$ such 
that the following diagram commute:
$$\xymatrix{A\ar[r]^{\pi}\ar[dr]_{f} & A/_mS\ar[d]^{!\tilde f} \\ & B}$$
where $\pi:A\rightarrow A/_mS$ is the canonical projection $\pi(a)=\overline a$.
 \end{prop}
 \begin{proof}
  Straightforward.
 \end{proof}
 
 \begin{prop}
  Let $A$ be a multiring, $I\subseteq A$ an ideal and $S\subseteq A$ a multiplicative subset such that $I\subseteq S$. Define 
$S/I=\{\overline s:s\in S\}$ (modulo $I$). Then $$(A/I)/_m(S/I)\cong A/_mS.$$
 \end{prop}
\begin{proof}
 Define $\Phi:A/I\rightarrow A/_mS$ given by $\Phi(\overline a^I)=\overline a^S$ and use the previous proposition. 
\end{proof}

\begin{prop}
 Let $A$ be a multiring and $P,S\subseteq A$ multiplicative subsets of $A$ such that $P\subseteq S$. Then
$$A/_mS\cong P^{-1}A/_mP^{-1}S.$$
\end{prop}
\begin{proof}
Straightforward.
\end{proof}

$ $

We discuss now the inductive limit of multirings. We will define this structure in terms of multivalued operations. However, is 
the case that this definition coincides with the categorical (and logical) definition of limit in terms of ternary relation.

Let $(I,\preceq)$ be a poset, and $\langle A_i,\varphi_{ij},I\rangle$ be a directed system, i.e, if $i,j\in I$ with $i\preceq j$, 
then there exist a morphism $\varphi_{ij}:A_i\rightarrow A_j$, $\varphi_{ii}=id_{A_i}$ and for $i\preceq j\preceq k$, the 
following diagram commutes:
$$\xymatrix{A_i\ar[dr]_{\varphi_{ij}}\ar[rr]^{\varphi_{ik}} & & A_k \\ & A_j\ar[ur]_{\varphi_{jk}} &}$$

Now, let $A$ be a multiring and $\langle A_i,\varphi_{ij},I\rangle$ be a directed system. A collection of morphisms 
$\{\psi_{i}:A_i\rightarrow A\}_{i\in I}$ is said to be compatible, if $\psi_j\varphi_{ij}=\psi_i$ for all $i\preceq j$.

\begin{prop}\label{multilimit}
 Let $\langle A_i,\varphi_{ij},I\rangle$ be a directed system of multirings. Then there exist a unique directed limit 
$\varinjlim\limits_{i\in I}A_i$ of the system.
\end{prop}
\begin{proof}
 The uniqueness is immediate. To show the existence, let 
 $$\varinjlim\limits_{i\in I}A_i=\dfrac{\dot\bigcup_{i\in I}A_i}{\sim},$$
 where $\sim$ is the following equivalence relation: for $x_i\in A_i$ and $x_j\in A_j$, $x_i\sim x_j$ if and only if there exist 
$k\succeq i,j$ such that $\varphi_{ik}(x_i)=\varphi_{jk}(x_j)$. We have $\varinjlim_{i\in I}A_i$ is a multiring with the obvious 
product and the sum defined by $\overline x_i+\overline y_j=\{\overline d:d\in\varphi_{ik}(x_i)+\varphi_{jk}(x_j),\,k\succeq 
i,j\}$. Moreover, $\psi:A_i\rightarrow\varinjlim_{i\in I}A_i$ given $\psi_i(x_i)=\overline x_i$ provides a compatible collection 
of morphisms. Then $\langle\varinjlim_{i\in I}A_i,\psi_i\rangle$ is the compatible system desired.
\end{proof}

 \begin{prop}
 $ $
  \begin{enumerate}[a -]
   \item  Let $\langle A_i,\varphi_{ij},I\rangle$ be a directed system of multirings and let $S_i\subseteq A_i$ be a multiplicative subset of 
$A_i$ for every $i\in 
I$. Suppose that for each $i, j \in I$ such that $i \preceq j$, $\varphi_{ij}[S_i] \subseteq S_j$. Then  $\langle 
S_i,{\varphi_{ij}}{\restriction},I\rangle$ and  $\langle A_i/_mS_i,\overline{\varphi}_{ij},I\rangle$ are directed systems, $\varinjlim_{i \in I} 
S_i$ is a multiplicative subset of $\varinjlim_{i \in I} A_i$ and
$$\varinjlim_{i \in I} A_i/_m S_i \cong (\varinjlim_{i \in I} A_i)/_m(\varinjlim_{i \in I} S_i).$$

   \item Let $\{A_i\}_{i\in I}$ be a family of multirings and $S_i\subseteq A_i$ be a multiplicative subset of $A_i$ for every $i\in 
I$. Then $$\prod_{i\in I}A_i/_m S_i\cong\prod_{i\in I}A_i/_m\prod_{i\in I} S_i.$$
  \end{enumerate}
 \end{prop}
 \begin{proof} Straitforward.
 \end{proof}

\subsection{Ordering Structures}
\hspace*{\parindent}
 
Part of the standard Artin-Schreier theory for fields can be extended to the multifield theory as in section 3 of \cite{marshall2006real}. For 
the convenience of the reader, we will list some results that we will use in the next sections:

\begin{defn}\label{defn:mfordering}
 Let $F$ be a multifield. A subset $P$ of $F$ is called an \textbf{ordering} if $P+P=\subseteq P$, $P\cdot P\subseteq P$, $P\cup 
-P=F$ and $P\cap -P=\{0\}$. The \textbf{real spectrum} of a multifield $F$, denoted $\mbox{Sper}(F)$, is defined to be the set of 
all orderings of $F$.
\end{defn}

\begin{defn}\label{defn:mfpreordering}
 A \textbf{preordering} of a multifield $F$ is defined to be a subset $T$ of $F$ satisfying $T+T\subseteq T$, $T\cdot T\subseteq 
T$ and $F^2\subseteq T$. Here, $F^2:=\{a^2:a\in F\}$. A multifield $F$ is said to be \textbf{real} if $-1\notin\sum F^2$. If $F$ 
is real, then $-1\ne1$. A preordering $T$ of $F$ is said to be \textbf{proper} if $-1\notin T$.
\end{defn}

\begin{prop}\label{prop:3.4marshall}
 Let $F$ be an multifield and $T$ a proper preordering of $F$. Then $T=\bigcap_{P\in X_T}P$, 
where $X_T=\{P\in\mbox{Sper}(F):T\subseteq P\}$.
\end{prop}
\begin{proof}
 Proposition 3.4 of \cite{marshall2006real}.
\end{proof}

$ $

Consider the multifield $Q_2$. The set $\{0,1\}$ is the unique ordering on $Q_2$. For any ordering $P$ on a multifield $F$, considering $\dot 
P:=P\setminus\{0\}$, we have $Q_P(F)=F/_m \dot P\cong Q_2$ by a unique isomorphism. So orderings of a multifield $F$ correspond bijectively to a 
multiring homomorphism 
$\sigma:F\rightarrow Q_2$ via $P=\sigma^{-1}(\{0,1\})$. We define $Q_{red}(F) = F /_m (\sum F^2\setminus\{0\})$.

\begin{prop}\label{prop:4.1marshall}
 For a real multifield $F$ are equivalent:
 \begin{enumerate}[a -]
  \item The multiring morphism $F\rightarrow Q_{red}(F)$ is an isomorphism;
  \item $\sum F^2=\{0,1\}$;
  \item For all $a\in F$, $a^3=a$ and $(a\in1+1)\Rightarrow(a=1)$.
 \end{enumerate}
\end{prop}
\begin{proof}
 Proposition 4.1 of \cite{marshall2006real}.
\end{proof}

\begin{defn}\label{defn:mfrealreduced}
 A multifield $F$ is said to be \textbf{real reduced} if satisfies the equivalent conditions of proposition 
\ref{prop:4.1marshall}.

A morphism of real reduced multifield is just a morphism of multifields. The category of real reduced multifields will be denoted 
by $\mathcal{MF}_{red}$.
\end{defn}

\begin{cor}\label{cor:4.2marshall}
 A multifield $F$ is real reduced if and only if $a^3=a$ for all $a\in F$ and $a\in1+1\Rightarrow a=1$.
\end{cor}

In general, for every multifield, the map $\mbox{sgn}:F\rightarrow Q_2^{\mbox{Sper}(F)}$ is available, and is given by the rule 
$\mbox{sgn}(x)=(\mbox{sgn}_P(x))_{P\in\mbox{Sper}(F)}$, where $$\mbox{sgn}_P(x)=\begin{cases}1\mbox{ if }x\in P\setminus(-P) \\ 
0\mbox{ if }x\in\mbox{supp}(P):=P\cap-P \\ -1\mbox{ if }x\in (-P)\setminus P\end{cases}$$

\begin{teo}[Local-Global principle]\label{cor:4.4marshall}
 For any real reduced multifield $F$, the map $\mbox{sgn}:F\hookrightarrow Q_2^{\mbox{Sper}(F)}$ is a strong embedding. In 
particular, for any real reduced multifield, $\mbox{Sper}(F)$ separate points of $F$ and $c\in a+b\subseteq F$ if and only if, 
for every $\sigma:F\rightarrow Q_2$, $\sigma(c)\in\sigma(a)+\sigma(b)$.
\end{teo}
\begin{proof}
 Proposition 4.4 \cite{marshall2006real}.
\end{proof}

\begin{defn}\label{defn:mrordering}
 Let $A$ be a multiring. A subset $P$ of $A$ is called an \textbf{ordering} if $P+P\subseteq P$, $P\cdot P\subseteq P$, $P\cup 
-P=A$ and $P\cap -P$ is a prime ideal of $A$, called the \textbf{support} of $P$. A set $T \subseteq A$ with just $T + T 
\subseteq T$, $T\cdot T \subseteq T$ and $A^2 \subseteq T$ is called a \textbf{preorder}. The \textbf{real spectrum} of a 
multiring $A$, denoted $\mbox{Sper}(A)$, is defined to be the set of all orderings of $A$.
\end{defn}

Like in the multifield case, orderings of a multiring $A$ correspond bijectively to multiring homomorphisms $\sigma:A\rightarrow 
Q_2$ via $P=\sigma^{-1}(\{0,1\})$. A \textbf{preordering} of a multiring $A$ is a subset $T$ of $A$ satisfying $T+T\subseteq T$, 
$TT\subseteq T$ and $A^2\subseteq T$. A preordering $T$ of $A$ is said to be \textbf{proper} if $-1\notin T$. Every ordering is a 
proper preordering. $\sum A^2$ us a preordering, and is the unique smallest preordering of $A$. A multiring $A$ is said to be 
\textbf{semireal} if $-1\notin\sum A^2$.

Let $A$ be a semireal multiring and consider $T$ a preorder. We denote by $Q_T (A)$ the image of $A$ in $Q_2 ^ 
{X_T}$, where $X_T = \{P \in \mbox{Sper}(A): T \subseteq P\}$.

Denote the image of $A$ in $Q_2^{X_T}$ by $Q_T(A)$. Addition on $Q_T(A)$ is defined by $\overline a+\overline b:=\{\overline 
c:c\in a+b\}$, $\overline a\overline b:=\overline{ab}$, $-\overline{a}:=\overline{-a}$. The zero element of $Q_T(A)$ is 
$\overline{0}$.

\begin{prop}[Local-Global principle]\label{prop:7.3marshall}
 Let $A$ be a semireal multiring and $T$ a proper preordering of $A$. Then:
 \begin{enumerate}
  \item $Q_T(A)$ is a multiring.
  \item $Q_T(A) \subseteq Q_2^{X_T}$ is a submultiring.
 \end{enumerate}
\end{prop}
\begin{proof}
 Proposition 7.3 \cite{marshall2006real}.
\end{proof}

We denote $Q_{\sum A^2}(A)$ by $Q_{red}(A)$ which we refer to as the \textbf{real reduced multiring} associated to $A$. 
Note that in multifield case, this reduction definition is compatible with the reduction for multifields, by theorem 
\ref{cor:4.4marshall}.

\begin{prop}\label{prop:7.5marshall}
 For a semireal multiring $A$, the map $a\mapsto\overline a$ from $A$ onto $Q_{red}(A)$ 
is an isomorphism if and only if $A$ satisfies the following properties:
\begin{enumerate}[a -]
 \item $a^3=a$.
 \item $a+ab^2=\{a\}$.
 \item $a^2+b^2$ contains a unique element.
\end{enumerate}
\end{prop}
\begin{proof}
 Proposition 7.5 \cite{marshall2006real}.
\end{proof}

\begin{defn}\label{defn:mrrealreduced}
 A semireal multiring satisfying one of the equivalent conditions of proposition \ref{prop:7.5marshall} will be called 
\textbf{real reduced multiring}. A morphism of real reduced multirings is just a morphism of multirings. The category of real 
reduced multirings will be denoted by $\mathcal{MR}_{red}$.
\end{defn}

\begin{cor}\label{cor:7.6marshall}
 A multiring $A$ is real reduced if and only if the following properties holds for all $a,b,c,d\in A$:
\begin{enumerate}[i -]
 \item $1\ne0$;
 \item $a^3=a$;
 \item $c\in a+ab^2\Rightarrow c=a$;
 \item $c\in a^2+b^2$ and $d\in a^2+b^2$ implies $c=d$.
\end{enumerate}
\end{cor}
\begin{proof}
 Corolllary 7.6 \cite{marshall2006real}.
\end{proof}

This implies that the morphism $a\mapsto\overline a$ from $A$ to $Q_{\mbox{red}}(A)$ is an isomorphism. In particular, follow by 
the local-global principle \ref{prop:7.3marshall} that for any real reduced multiring $A$, $c\in a+b\subseteq F$ if and only if, 
for every $\sigma:A\rightarrow Q_2$, $\sigma(c)\in\sigma(a)+\sigma(b)$.

\section{Multifields, Abstract Ordering Spaces, and Special Groups}

Marshall's abstract space of orderings (AOS), arised in the decade of 80's and presented in \cite{marshall1996spaces}, are 
important because they generalize both theory of orderings on fields and the reduced theory of quadratic forms. But only in the 
decade of 90s a (finitary) first-order the theory has arisen that generalizes the reduced and non-reduced theory of quadratic 
forms simultaneously. This is the theory of Special Groups of M. Dickmann and F. Miraglia, presented in 
\cite{dickmann2000special} (and generalized for invertible coefficients in rings in \cite{dickmann2015faithfully}). 

The reduced special groups are dually equivalent to the abstract ordering spaces ($\mathcal{RSG}\simeq\mathcal{AOS}^{op}$). This 
simplicity brings two new methods and tools to the algebraic theory of quadratic forms, the $K$-theory 
(\cite{dickmann2006algebraic}) and the boolean hull of a special groups (\cite{dickmann2000special}), culminating in 
a proof of problems on quadratic forms theories that were open by 25 years (Marshall's and Lam signature conjectures, 
\cite{dickmann1998quadratic} and \cite{dickmann2003lam} respectively).

%

From the middle of years 2000's, the development of special group theory focused on expanding the class of examples (see 
for example, \cite{dickmann2008quadratic} and \cite{dickmann2009representation}). In recent book of Dickmann and Miraglia 
\cite{dickmann2015faithfully}, they extend the classical algebraic theory of quadratic forms over fields to a broad class of 
commutative rings with unit (of course, 
which was mediated by the theory of special groups). The context is of a ring $A$ of characteristic not 2, with $-1\notin\sum A^2$ 
and $2\in\dot A$.

Given a such ring $A$ and a preordering $T$ on $A$, they define that two $n$-dimensional forms $\varphi=a_1X_1^2+...+a_nX_n^2$, 
$\psi=b_1X_1^2+...+b_nX_n^2$ with $a_i,b_i\in\dot A$ are $T$-isometric\index{isometric}, $\varphi\approx_T\varphi$ if there is a 
sequence $\varphi_0,\varphi_1,...,\varphi_k$ of $n$-dimensional diagonal forms over $\dot A$, such that $\varphi=\varphi_0$, 
$\psi=\varphi_k$ and for every $1\le i\le k$, $\varphi_i$ is either isometric to $\varphi_{i-1}$ in the usual sense that there is 
a matrix $M\in\mbox{GL}_n(A)$ such that $\varphi_i=M\varphi_{i-1}M^t$ or there are $t_1,...,t_n\in\dot T$ such that 
$\varphi_i=\langle t_1x_1,...,t_nx_n\rangle$ and $\varphi_{i-1}=\langle x_1,...,x_n\rangle$. Value representation relation $D_T$ 
on $(A,T)$ is given by: for $a,b_1,...,b_n\in\dot A$,
$$a\in D^v_T(b_1,...,b_n)\Leftrightarrow\exists\,t_1,...,t_n\in T\mbox{ such that }a=\sum\limits^n_{i=1}t_ib_i.$$

Given a preordered ring $(A,T)$, they associate a structure $G_T(A)$, whose domain is $\dot A/\dot T$, endowed with the product 
operation induced by $\dot A$, togheter with a binary isometry relation $\equiv_{G_T(A)}$, defined on ordered pairs of elements 
of $\dot A/\dot T$, and having $-1=-1/\dot T$ as distinguished element. The structure $(G_T(A),\equiv_{G_T(A)},-1)$ is not quite 
a special group, but satisfy SG0, SG1, SG2, SG3 and SG5. They observed that the ring-theoretic approach, based on the definition 
of $n$-isometry and the formal approach via $G_T(A)$, though related, are far from identical.

Beside this, they called \textbf{$T$-faithfully quadratic} any preordered ring $(A,T)$ such that $G_T(A)$ is a special group and 
$T$-isometry and value representation in $(A,T)$ are faithfully coded by the corresponding formal notions in $G_T(A)$. After this 
brilliant idea, they was able to replicate most of the consequences of the theory of special groups in field theory in this 
extended ring-theoretic context.

In this section, we will construct two equivalence of categories. The first one, is given by a functor 
$\mathcal{MF}_{red}\rightarrow\mathcal{AOS}^{op}$ where $\mathcal{AOS}$ is the category of abstract ordering spaces as in 
\cite{marshall1996spaces} and $\mathcal{MF}_{red}$ is the category of real reduced multifield. This functor was indicated by 
Marshall in \cite{marshall2006real}, so we will establish precisely its definition and examine some properties. The second 
equivalence extends the first one and is given by a functor $\mathcal{SG}\rightarrow\mathcal{SMF}$, from the category of special 
groups to the category of special multifields, here introduced.

Recall that for any set $X$, we give to $\{-1,1\}^X$ a group structure by defining $(ab)(x)=a(x)b(x)$, and that if $G$ is a 
group of exponent 2, the \textbf{character group} of $G$ is the group $\chi(G):=\mbox{Hom}(G,\mathbb Z_2)$.

\begin{defn}[Space of Orderings]\label{defn:aos}
 An \textbf{abstract ordering space or space of orderings}, is a pair $(X,G)$ satisfying:
\begin{description}
 \item [AX1 -] $X$ is a non-empty set, $G$ is a subgroup of $\{-1,1\}^X$, $G$ contains the constant function $-1$, and $G$ 
separates points in $X$ (i.e, if $x,y\in X$, $x\ne y$, then there exists $a\in G$ such that $a(x)\ne a(y)$).
\end{description}

If $a,b\in G$ we define the {\bf value set} $D(a,b)$ to be the set of all $c\in G$ such that for each $x\in X$ either $c(x)=a(x)$ 
or $c(x)=b(x)$. In particular, $a$ and $b$ are both elements of $D(a,b)$.

\begin{description}
\item [AX2 -] If $x\in\chi(G)$ satisfies $x(-1)=-1$ and $a,b\in\mbox{ker}(x)\Rightarrow D(a,b)\subseteq\mbox{ker}(x)$, then $x$ 
is in the image of the natural embedding $X\hookrightarrow\chi(G)$.

\item [AX3 (Associativity) -] For all $a,b,c\in G$, if $t\in D(a,r)$ for some $r\in D(b,c)$ then $t\in D(s,c)$ for some $s\in 
D(a,b)$.
\end{description}
\end{defn}

\begin{defn}\label{defn:aormorphism}
 A \textbf{morphism} of abstract ordering spaces $\alpha:(X,G)\rightarrow(Y,H)$  
is a mapping $\alpha:X\rightarrow Y$ such that for each $h\in H$, the composite function $h\circ\alpha:X\rightarrow\{-1,1\}$ is 
an element of $G$. Note that this implies that $\alpha$ induces a group homomorphism 
$h\mapsto h\circ\alpha$ from $H$ to $G$. Also, $\alpha^{-1}(U(h))=U(h\circ\alpha)$ for each $h\in H$, so $\alpha$ is continuous.

An isomorphism from $(X, G)$ to $(Y, H)$ is an AOS-morphism $\alpha:X\rightarrow Y$ which is bijective and such that the induced 
group homomorphism $h\mapsto h\circ\alpha$ is also bijective.
\end{defn}

\begin{teo}\label{teo:aortomfred}
 Let $(X,G)$ be a space of orderings and set $M(G)=G\cup\{0\}$ where $0:=\{G\}$. Then $(M(G),+,\cdot,-,0,1)$ is a real reduced 
multifield with the extended operations: 
\begin{itemize}
   \item $a\cdot b=\begin{cases}0\,\mbox{if }a=0\mbox{ or }b=0 \\ a\cdot 
b\,\mbox{otherwise}\end{cases}$
   \item $-(a)=(-1)\cdot a$
   \item $a+b=\begin{cases}\{b\}\,\mbox{if }a=0 \\ \{a\}\,\mbox{if }b=0 \\ M(G)\,\mbox{if 
}a=-b,\,\mbox{and }a\ne0 
\\ 
D(a,b)\,\mbox{otherwise}\end{cases}$
  \end{itemize}
\end{teo}
\begin{proof}
 Firstly, observe that $+$ is well-defined. Then, we will verify the conditions of definition \ref{defn:multiring}:
\begin{enumerate}[i -]
 \item For this, we will check the conditions of definition \ref{defn:multigroupI}.  
 \begin{enumerate}[a -]
  \item We will prove that if $d \in a + b$, then $a \in d + (-b)$ and $b \in (-a) + d$.
  If $a=0$ or $a=-b$, then $d\in a+b$ implies trivially that $a\in d+(-b)$ and $b\in(-a)+d$. Now, let $a,b\ne0$ with 
$a\ne-b$ (this implies $d\ne0$). Now it is enough to note that for all $x \in X$ if $d(x)\in\{a(x),b(x)\}$, then $a(x)\in\{d(x),-b(x)\}$ and 
$b(x) \in \{-a(x),d(x)\}$.
Thus $a \in d + (-b)$ and $b \in (-a) + d$.

  \item $(y\in x+0)\Leftrightarrow(x=y)$ is direct consequence of the definition of multivaluated sum.
   
   \item $a+0 = \{a\} = 0+a$ and $a+(-a)=M(G)=(-a)+a$. Let $a,b\in M(G)$, $a,b\ne0$ and $a\ne -b$. Since $D(a,b)=D(b,a)$, we have $a+b=b+a$. 
Then, the commutativity holds.

  \item Now we prove the associativity, that is, $(a + b) + c = a + (b + c)$ for all $a,b,c \in M(G)$. If $a=0$ (the cases $b=0$ and $c=0$ are 
analogous), then $0+(b+c)=\{0+g:g\in b+c\}=b+c$ 
and $(0+b)+c=(\{b\})+c=b+c$.

Now, assume that $a,b,c\ne0$. If $-c \in a + b$, then

$$(a+b)+c=\bigcup\{g+c:g\in a+b\}=M(G)\,\mbox{(I)}$$

and

$$a+(b+c)=\bigcup\{a+h:h\in b+c\}=M(G)\,\mbox{(II)}$$

because $-a\in b+c$. So $\mbox{(I)}=\mbox{(II)}$. For the case 
$a,b,c\ne0$, $-c \notin a + b$ we have

\begin{align*}
 (a+b)+c & = \bigcup\{g + c : g \in a + b\} \\
	 & = \bigcup\{g + c : g \in D(a,b)\} \ = \bigcup_{g \in D(a,b)} D(c,g) \, \mbox{(III)}
\end{align*}

because the hypothesis implies $a \neq -b$ and

\begin{align*}
 a+(b+c) & = \bigcup\{a + h : h \in b + c\} \\ 
	 & = \bigcup\{a + h : h \in D(b,c)\} \ = \bigcup_{h \in D(b,c)} D(a,h)\,\mbox{(IV)}
\end{align*}

because by item a) above $-a \notin b + c$ and so $b \neq -c$.
By the inductive description of the value sets (as in 2.2 of \cite{marshall1996spaces}) we have $\mbox{(III)}=\mbox{(IV)}$. Thus
$(a+b)+c=a+(b+c)$ for all $a,b,c\in M(G)$.
 \end{enumerate}
 
   \item Since $(G,\cdot,1)$ is an abelian group, we conclude that $(M(G),\cdot,1)$ is a commutative monoid. Beyond this, every nonzero element 
of $M(G)$ has an inverse.
  
  \item $a\cdot 0=0$ for all $a\in M(G)$ is direct from definition.
  
  \item For the distributive property, let $a,b,c,d \in M(G)$ with $d \in a + b$. If $a = 0$ or $b = 0$ or $c = 0$, then $cd \in ca + cb$ (*) and 
if $a,b \neq 0$ and $a = - b$ we also
  conclude (*) easily. If $a,b \neq 0$ and $a \neq -b$, then by definition $d \neq 0$ and for all $x \in X$ we have $d(x) \in \{a(x), b(x)\}$. 
Then $d(x)c(x) \in \{a(x)c(x), b(x)c(x)\}$
  and so $cd \in D(ca,cb)$. Therefore $cd \in ca + cb$.
\end{enumerate}

Then, $(M(G),+,-,\cdot,0,1)$ is a multifield. As $G$ is a subgroup of $\{-1,1\}^X$, we have $G$ is a group of exponent 2, i.e, 
$g^2=1$ for all $g\in G$ and then, $a^3=a$ for all $a\in M(G)$. If $a\in1+1$, then $a(x)=1$ for all $x\in X$. This implies $a=1$ and so 
$(M(G),+,-,\cdot,0,1)$ is a real reduced multifield.
\end{proof}

\begin{cor}\label{cor:equiv3}
 The correspondence $(X,G) \mapsto M(G)$ extends to a contravariant functor $M:\mathcal{AOS}^{op}\rightarrow\mathcal{MF}_{\mbox{red}}$.
\end{cor}

\begin{proof}
Let $(X,G)$ and $(Y,H)$ abstract ordering spaces and $\alpha:Y\rightarrow X$ be an AOS-morphism. By definition 
\ref{defn:aormorphism}, $\alpha$ induces a group homomorphism $\varphi:G\rightarrow H$ given by $\varphi(g)=g\circ\alpha$. Define 
$M(\alpha)=\tilde{\varphi}:M(G)\rightarrow M(H)$ extending this morphism $\varphi$ to $M(G)$ making $\tilde{\varphi}(0)=0$. Note 
that we already have $\tilde{\varphi}(1)=1$, $\tilde{\varphi}(-1)=-1$ and $\tilde{\varphi}(ab) = \tilde{\varphi}(a) \tilde{\varphi}(b)$ for all 
$a,b \in M(G)$.

Then, we just need to prove that for all $a,b,c\in M(G)$, $c \in a+b \Rightarrow \tilde{\varphi}(a) \in \tilde{\varphi}(b) + \tilde{\varphi}(c)$. 
The cases $a = 0, b=0, c=0$ and 
$(a,b \neq 0, \varphi(a) = -\varphi(b))$ are trivial and we can assume $a,b,c \neq 0, \varphi(a) \neq -\varphi(b)$ (in particular, $a \neq -b$). 
Hence, we need to prove that $c\in 
D(a,b)\Rightarrow\varphi(c)\in D(\varphi(a),\varphi(b))$. But note that

\begin{align*}
 c\in D(a,b)\Rightarrow c(x)=a(x)\vee c(x)=b(x)\,\forall\,x\in X\Rightarrow \\
 c(\alpha(y))=a(\alpha(y))\vee c(\alpha(y))=b(\alpha(y))\,\forall\,y\in Y\Rightarrow \\
 c\circ\alpha\in D(a\circ\alpha,b\circ\alpha)\Rightarrow D(\varphi(a),\varphi(b)).
\end{align*}

Therefore $M(\varphi)$ is a MF-morphism. If $\xymatrix{(Z,K)\ar[r]^{\beta} & (Y,H)\ar[r]^{\alpha} & (X,G)}$ are AOS-morphism, 
with $\varphi:G\rightarrow H$ and $\tau:H\rightarrow K$ the respectively induced group homomorphisms, the fact of 
$M(\alpha\beta)=M(\beta)M(\alpha)$ is direct consequence of $\alpha\beta$ be an AOS-morphism.
\end{proof}

$ $

Let $F$ be an real reduced multifield. Observe that by the local-global principle for multifield 
\ref{cor:4.4marshall} we have the following identities:
\begin{itemize}
 \item $a\in a+b$;
 \item If $a\ne0$, then $a+(-a)=F$;
 \item $a\ne0\Rightarrow\sigma(a)\ne0$ for all $\sigma\in\mbox{Sper}(F)$.
\end{itemize}

Consider the map $sgn \colon F \to Q_2^{\mbox{sper}(F)}$ from $\ref{cor:4.4marshall}$.

\begin{teo} \label{mfAos}
 Let $F$ be a real reduced multifield. Then $(\mbox{sper}(F), sgn(\dot{F}))$ is an abstract order space.
\end{teo}

\begin{proof}
 Since $sgn \colon F \to Q_2^{\mbox{sper}(F)}$ is a injective morphism, it is easy to see that $sgn(\dot{F}) \subseteq \{1,-1\}^{\mbox{sper}(F)}$ 
is subgroup that contains
 the constant functions and $\mbox{sper}(F)$ separates points of $sgn(\dot{F})$.
 
 \begin{claim*}
  Let $a,b,c \in \dot{F}$ and $a' = sgn(a), b' = sgn(b), c' = sgn(c)$. Then $a' \in D(b',c')$ if, and only if, $a \in b + c$.
 \end{claim*}

 \begin{proof}
 By definition, we have
 
 \begin{align*}
  a' \in D(b',c') & \Leftrightarrow \mbox{for all } \sigma \in \mbox{sper}(F) \colon a'(\sigma) \in \{b'(\sigma), c'(\sigma)\} \\
		  & \Leftrightarrow \mbox{for all } \sigma \in \mbox{sper}(F) \colon \sigma(a) \in \{\sigma(b), \sigma(c) \} \\
		  & \Leftrightarrow \mbox{for all }\sigma \in \mbox{sper}(F) \colon \sigma(a) \in D_{Q_2} (\sigma(b), \sigma(c)) \\
		  & \Leftrightarrow a \in b + c,
 \end{align*}

 where the last equivalence follows by the local-global principle ($\ref{cor:4.4marshall}$).
 \end{proof}
 
 Now we check the axioms $AX2, AX3$ from definiton $\ref{defn:aos}$.
 
 \begin{description}
  \item [AX2] Let $x \in \chi(sgn(\dot{F}))$ satisfying $x(-1) = -1$ and for all $a',b' \in sgn(\dot{F})$ if $a',b' \in ker(x)$, then $D(a',b') 
\subseteq ker(x)$. 
  Using the above claim, since $sgn$ is a strong injective morphism, exists $\tilde{\sigma} \colon \dot{F} \to \{1,-1\}$ satisfying
  
  \begin{enumerate}[i)]
   \item $\tilde{\sigma}(-1) = -1$ and $\tilde{\sigma}(ab) = \tilde{\sigma}(a)\tilde{\sigma}(b)$ for all $a,b \in \dot{F}$.
   
   \item If $a,b \in ker(\tilde{\sigma})$, then $\dot{(a + b)} \subseteq ker(\tilde{\sigma})$.
   
   \item $\tilde{\sigma} = x \circ sgn$.
  \end{enumerate}

  If we define $\sigma \colon F \to Q_2$ by $\sigma(a) = \tilde{\sigma} (a)$ for $a \neq 0$ and $\sigma(0) = 0$, it is easily verifiable that 
$\sigma$ satisfies every axiom
  of multiring morphism except possibly the multivalued sum preservation. Thus let $a,b,c \in F$ with $a \in b + c$. We want to prove that 
$\sigma(a) \in \sigma(b) + \sigma(c)$ (*).
  The cases $0 \in \{a,b,c\}$ or $\sigma(b) = - \sigma(c)$ are trivial . So assume $a,b,c \neq 0$ and $\sigma(b) = \sigma(c)$. If $\sigma(b) = 
\sigma(c) = 1$, then by 
  item $ii)$ above $\sigma(a) = 1$ and so (*) is verified. If $\sigma(b) = \sigma(c) = -1$, then since $-a \in -b -c$, we conclude by the last 
case that $\sigma(-a) = 1$, that is, 
  $\sigma(a) = -1$. Thus (*) is fully verified and $\sigma \in \mbox{sper}(F)$. Therefore by the item $iii)$ above, $x \in \chi(sgn(\dot{F}))$ is 
the image of $\sigma$ by the
  canonical map $X \to \chi(sgn(\dot{F}))$.
  
  \item [AX3] Now let $t' \in D(a',r')$ with $r' \in D(b',c')$. We want to prove that exist $s' \in D(a',b')$ such that $t' \in D(s',c')$.
  Let $t,a,r,b,c \in \dot{F}$ with $t' = sgn(t), a' = sgn(a), r' = sgn(r), b' = sgn(b), c' = sgn(c)$. By the above claim, $t \in a + r, r \in b + 
c$. Using the associative 
  property of $F$, exist $s \in a + b$ with $t \in s + c$. If $s \neq 0$, then again by the claim $sgn(s) \in D(a',b')$ and $t' \in D(sgn(s), 
c')$. If $s = 0$, then $t = c$
  and so choosing $s' = a'$ we have $s' \in D(a',b')$ and $t' \in D(s',c')$.
 \end{description}

\end{proof}

\begin{cor}
 The correspondence $F \to (\mbox{sper}(F), sgn(\dot{F}))$ extends to a contravariant functor $\mbox{sper} \colon \mathcal{MF}_{\mbox{red}} \to 
\mathcal{AOS}^{op}$.
\end{cor}

\begin{proof}
 Let $f \colon F \to G$ a morphism between real reduced mutifields. Consider the induced map $\alpha_f \colon\mbox{sper}(G) \to \mbox{sper}(F)$. 
Note that
 given $sgn_F(a) \in sgn_F(\dot{F})$, $sgn_F(a) \circ \alpha_f = sgn_G(f(a)) \in sgn_G(\dot{G})$ and so 
 $\alpha_f \colon (\mbox{sper}(G), sgn_G (\dot{G})) \to (\mbox{sper}(F), sgn_F (\dot{F}))$ is a AOS-morphism. Furthermore, it is easy to see that 
if
 $f \colon F \to G$ and $g \colon G \to H$ are morphism between real reduced multifields, then $\alpha_{g \circ f} = \alpha_f \circ \alpha_g$.
\end{proof}

\begin{teo}
 The functors $M:\mathcal{AOS}^{op}\rightarrow\mathcal{MF}_{\mbox{red}}$ and $\mbox{sper} \colon \mathcal{MF}_{\mbox{red}} \to 
\mathcal{AOS}^{op}$ establish a equivalence of categories.
\end{teo}

\begin{proof}
 Let $F$ a real reduced multifield and $(X,G)$ an abstract order space. We exhibit two ismorphisms: $\eta \colon F \to M(sgn(\dot{F}))$ 
multifield morphism and 
 $\alpha \colon (X,G) \to (\mbox{sper}(M(G)), sgn(\dot{M(G)}))$ AOS-morphism. 
 
 \begin{description}
  \item [$\eta$]: Consider the strong injective morphism $sgn \colon F \to Q_2^{\mbox{sper}(F)}$. By the claim in the theorem $\ref{mfAos}$, in 
the AOS 
  $(\mbox{sper}(F), sgn(\dot{F}))$ given $a,b \in \dot{F}$ the value set is given by $D(sgn(a), sgn(b)) = sgn(\dot{a + b})$. Define $\eta \colon 
F \to M(\dot{F})$ by

  \begin{equation*}
    \eta(a) = 
    \begin{cases}
      sgn(a) &\text{if $a \neq 0$}\\
      0      & 	\text{if $a = 0$}
    \end{cases}
  \end{equation*}
  
  It is imediate that $\eta$ is bijective and preserves constants and product. Note that given $a,b \in F$, we have by definition
  
  \begin{enumerate}
   \item [.] If $a = 0$, then $\eta(a + b) = \eta(b) = \eta(a) + \eta(b)$.
   
   \item [.] If $a,b \neq 0$ and $a = - b$, then $\eta(a + b) = \eta(F) = M(sgn(\dot{F})) = \eta(a) + \eta(b)$.
   
   \item [.] If $a,b \neq 0$ and $a \neq -b$, then 
   $$\eta(a + b) = sgn(a + b) = D(sgn(a), D(sgn(b))) = \eta(a) + \eta(b).$$
  \end{enumerate}

  Thus $\eta$ is an multifield isomorphism.
    
  \item [$\alpha$]: Given $x \in X$, we can consider the evaluation function $x \colon G \to \{1,-1\}$ given by $x(g) = g(x)$ for all $g \in G$. 
Then consider $\sigma_x \colon M(G) \to Q_2$
  given by
  
  \begin{equation*}
    \sigma_x (g) = 
    \begin{cases}
      g(x) &\text{if $g \neq 0$}\\
      0      & 	\text{if $g = 0$}
    \end{cases}
  \end{equation*}
  
  Following the proof of item $AX2$ in theorem $\ref{mfAos}$, we see that $\sigma_x$ is a multiring morphism. Define
  $\alpha \colon X \to \mbox{sper}(M(G))$ by $\alpha(x) = \sigma_x$. Note that given $sgn(g) \in sgn(\dot{M(G)})$, $sgn(g) \circ \alpha = g \in 
G$ (*) and so $\alpha$ is in fact 
  an AOS-morphism. The injectivity of $\alpha \colon X \to \mbox{sper}(M(G))$ is a consequence of $X$ separate points of X and the surjective 
follows by the axiom
  $AX2$ of definition $\ref{defn:aos}$. The equality (*) shows that $\alpha$
  is in fact an isomorphism, as desired.
  
  \end{description}

\end{proof}

Let $A$ be a set and $\equiv$ a binary relation on $A\times A$. We extend $\equiv$ to a binary relation $\equiv_n$ on $A^n$, by 
induction on $n\ge2$, as follows:
\begin{enumerate}[i -]
 \item $\equiv_2=\equiv$.
 \item $\langle a_1,...,a_n\rangle\equiv_n\langle b_1,...,b_n\rangle$ if and only there are 
 $x,y,z_3,...,z_n\in A$ such that $\langle a_1,x\rangle\equiv\langle b_1,y\rangle$,
 $\langle a_2,...,a_n\rangle\equiv_{n-1}\langle x,z_3,...,z_n\rangle$ and
 $\langle b_2,...,b_n\rangle\equiv_{n-1}\langle y,z_3,...,z_n\rangle$.
\end{enumerate}

Whenever clear from the context, we frequently abuse notation and indicate the afore-described extension 
$\equiv$ by the same symbol.  

\begin{defn}[Special Group]\label{defn:sg}
 A \textbf{special group} is an tuple $(G,-1,\equiv)$, where $G$ is a group of exponent 2, 
i.e, $g^2=1$ for all $g\in G$; $-1$ is a distinguished element of $G$, and $\equiv\subseteq G\times 
G\times G\times G$ is a relation (the special relation), satisfying the following axioms for all 
$a,b,c,d,x\in G$:
\begin{description}
 \item [SG 0 -] $\equiv$ is an equivalence relation on $G^2$;
 \item [SG 1 -] $\langle a,b\rangle\equiv \langle b,a\rangle$;
 \item [SG 2 -] $\langle a,-a\rangle\equiv\langle1,-1\rangle$;
 \item [SG 3 -] $\langle a,b\rangle\equiv\langle c,d\rangle\Rightarrow ab=cd$;
 \item [SG 4 -] $\langle a,b\rangle\equiv\langle c,d\rangle\Rightarrow\langle 
a,-c\rangle\equiv\langle-b,d\rangle$;
 \item [SG 5 -] 
$\langle a,b\rangle\equiv\langle c,d\rangle\Rightarrow\langle ga,gb\rangle\equiv\langle 
gc,gd\rangle,\,\mbox{for all }g\in G$.
 \item [SG 6 (3-transitivity) -] the extension of $\equiv$ for a binary relation on $G^3$ is a 
transitive relation.
\end{description}
\end{defn}

A group of exponent 2 satisfying SG0-SG5 is called \textbf{pre-special group}. A pre-special group (or special group) 
$(G,-1,\equiv)$ is \textbf{reduced} if $1\ne-1$ and if $\langle a,a\rangle\equiv\langle1,1\rangle\Rightarrow a=1$.

A \textbf{$n$-form} (or form of dimension $n\ge1$) is an $n$-tuple of elements of $G$. An 
element $b\in G$ is \textbf{represented} on $G$ by the form $\varphi=\langle a_1,...,a_n\rangle$, 
in symbols $b\in D_G(\varphi)$, if there exists $b_2,...,b_n\in G$ such that $\langle 
b,b_2,...,b_n\rangle\equiv\varphi$. Now, some examples:

\begin{ex}[The trivial special relation]\label{ex2.2}
 Let $G$ be a group of exponent 2 and take $-1$ as any element of $G$ different of 1. For $a,b,c,d\in G$, 
define $\langle a,b\rangle\equiv_t\langle c,d\rangle$ if and only if $ab=cd$. Then $G_t=(G,\equiv_t,-1)$ 
is a SG (\cite{dickmann2000special}).
\end{ex}

\begin{ex}[Special group of a field]\label{ex2.3}
 Let $F$ be a field. We denote $\dot F=F\setminus\{0\}$, $\dot F^2=\{x^2:x\in 
\dot F\}$ and $\Sigma \dot F^2=\{\sum_{i\in I}x_i^2:I\mbox{ is finite and }x_i\in 
\dot F^2\}$. Let $G(F)=\dot F/\dot F^2$. In the case of $F$ is be formally real, we have
$\Sigma \dot F^2$ is a subgroup of $\dot F$, then we take $G_{\mbox{red}}(F)=\dot F/\Sigma 
\dot F^2$. Note that $G(F)$ and $G_{\mbox{red}}(F)$ are groups of exponent 2. In 
\cite{dickmann2000special} they prove that $G(F)$ and $G_{\mbox{red}}(F)$ are special groups with the special 
relation given by usual notion of isometry, and $G_{\mbox{red}}(F)$ is always reduced.
\end{ex}

\begin{defn}\label{defnmorph}
 A map $\xymatrix{(G,\equiv_G,-1)\ar[r]^f & (H,\equiv_H,-1)}$ between pre-special groups is a \textbf{morphism 
of pre-special groups or PSG-morphism} if $f:G\rightarrow H$ is a homomorphism of groups, $f(-1)=-1$ and for all 
$a,b,c,d\in G$
$$\langle a,b\rangle\equiv_G\langle c,d\rangle\Leftrightarrow
\langle f(a),f(b)\rangle\equiv_H\langle f(c),f(d)\rangle$$
A \textbf{morphism of special groups or SG-morphism} is a PSG-morphism between the correspondents pre-special groups. $f$ 
will be an isomorphism if is bijective and $f,f^{-1}$ are PSG-morphisms. 
\end{defn}

The category of special groups (respectively reduced special groups) and theirs morphisms will be denoted by $\mathcal{SG}$ 
(respectively $\mathcal{RSG}$). Now, we will analyze the connections between the $\mathcal{SG}$ and $\mathcal{MF}$. For this, we 
need more results about special groups and their characterization. For this, we use the results proved in Lira's thesis 
\cite{de1996groupes}. Consider these axioms concerns about a group of exponent 2 with a distinguished element: 

\begin{description}
 \item [SG 7 -] $\forall a\,\forall a'\,\forall x\,\forall t\,\forall t'\,\forall 
y[(a,a')\equiv(x,t)\wedge(t,t')\equiv(1,y)]$ \newline $\Rightarrow\exists a''\,\exists s\,\exists 
s'[(a,a'')\equiv(y,s)\wedge(s,s')\equiv(1,x)]$.

An equivalent statement for SG7 is
$$\bigcup_{t\in D_G(1,y)}D_G(x,t)=\bigcup_{s\in D_G(1,x)}D_G(y,s)$$
for all $x,y\in G$.

 \item [SG 8 -] For all forms $f_1,...,f_n$ of dimension 3 and for all 
$a,a_2,a_3,b_2,b_3\in G$,
$$\langle a,a_2,a_3\rangle\equiv f_1\equiv...\equiv f_n\equiv\langle a,b_2,b_3\rangle\Rightarrow\langle 
a_2,a_3\rangle\equiv\langle b_2,b_3\rangle.$$
 \item [SG 9 -] $\forall a\,\forall b\,\forall c\,\forall 
d[\langle a,b,ab\rangle\equiv\langle c,d,cd\rangle\Rightarrow\langle a,b,ab\rangle\equiv\langle 
d,c,cd\rangle]$
\end{description}

\begin{prop}\label{sg6moreasy}
 Let $(G,-1, \equiv)$ be a pre-special group. The following are equivalent:
\begin{enumerate}[i -]
 \item $G \models SG6$
 \item $G \models SG7\wedge SG8$
 \item $G \models SG9$
\end{enumerate}
\end{prop}
\begin{proof}
 \cite{de1996groupes} page 32.
\end{proof}

\begin{prop}\label{sg.to.mf}
 Let $(G,\equiv,-1)$ be a special group and define $M(G)=G\cup\{0\}$ where $0:=\{G\}$\footnote{Here, 
the choice of the zero element was ad hoc. Indeed, we can define $0:=\{x\}$ for any $x\notin G$.}. Then 
$(M(G),+,-,\cdot,0,1)$ is a multifield, where 
\begin{itemize}
   \item $a\cdot b=\begin{cases}0\,\mbox{if }a=0\mbox{ or }b=0 \\ a\cdot 
b\,\mbox{otherwise}\end{cases}$
   \item $-(a)=(-1)\cdot a$
   \item $a+b=\begin{cases}\{b\}\,\mbox{if }a=0 \\ \{a\}\,\mbox{if }b=0\\ M(G)\,\mbox{if 
}a=-b,\,\mbox{and }a\ne0 
\\ 
D_G(a,b)\,\mbox{otherwise}\end{cases}$
  \end{itemize}
\end{prop}
\begin{proof}
 The proof is analogous to \ref{teo:aortomfred}.
\end{proof}

\begin{cor}\label{cor:equiv1}
 The correspondence $G\mapsto M(G)$ extends to a faithful functor $M:\mathcal{SG}\rightarrow\mathcal{MF}$. 
\end{cor}
\begin{proof}
 Let $f:G\rightarrow H$ be a SG-morphism. We will extend $f$ to $M(f):M(G)\rightarrow M(H)$ by 
$M(f)\upharpoonleft_G=f$ and $M(f)(0)=0$. By the definition of SG-morphism we have $M(f)(1)=1$, 
$M(f)(-a)=-a$ and $M(f)(ab)=M(f)(a)M(f)(b)$. As $d\in D_G(a,b)$ implies $f(d)\in D_H(f(a),f(b))$ we have 
$d\in a+b\Rightarrow M(f)(d)\in M(f)(a)+M(f)(b)$ for all $a,b\in M(G)$. So $M(f)$ is a multiring 
morphism. Now, let $\xymatrix{ G\ar[r]^{f} & H\ar[r]^{g} & K}$ be SG-morphisms. How $M(f\circ 
g)\upharpoonleft_G=f\circ g=M(f)\upharpoonleft_G\circ M(g)\upharpoonleft_G$ and $M(f\circ 
g)(0)=0=M(f)\circ M(g)(0)$, we have $M(f\circ g)=M(f)\circ M(g)$. Since $M(id)=id$, then 
$M:\mathcal{SG}\rightarrow\mathcal{MF}$ is a functor.

This functor is faithful, because if $G$ and $H$ are special groups and $f,g:G\rightarrow H$ are 
SG-morphisms such that $M(f),M(g):M(G)\rightarrow M(H)$ are equal, then 
$M(f)|_{M(G)\setminus\{0\}}=M(g)|_{M(G)\setminus\{0\}}$ and therefore $f=g$, since 
$M(G)\setminus\{0\}=G$.   
\end{proof}

\begin{prop}\label{prop:missues}
 Let $G$ be an SG and $M(G)$ as above. Then:
 \begin{enumerate}[i -]
  \item $a^2=1$ for all $a\in M(G)\setminus\{0\}$;
  \item $1 \in 1 + a$ for all $a\in M(G)$;
  \item $1 + a$ is closed by multiplication for all $a \in M(G)$;
  \item If exist $p \in \dot M(G)$ such that
  
  \begin{align*}
   a & \in c + cp \\
   b & \in p + ap \\
   d & \in p + cp.
  \end{align*}

  then exist $l \in\dot M(G)$ such that
  
  \begin{align*}
   a & \in d + dl \\
   b & \in l + al \\
   c & \in l + dl.
  \end{align*}
  
 \end{enumerate}
\end{prop}
\begin{proof}
 $ $
 \begin{enumerate}[i -]
  \item Is just the fact of $G$ be a group of exponent 2.
  \item Trivial.
  \item If $a = 0$ or $a = -1$ it is trivial. If $a \neq 0, -1$, given $x, y \in 1 + a = D_G(1,a)$, we have 
$\langle x, xa \rangle \equiv \langle 1,a \rangle $ and $\langle y, ya\rangle \equiv \langle 1, a\rangle$. 
Multiplying the first equality by one, we have $\langle xy, xya\rangle \equiv \langle y, ya\rangle \equiv 
\langle 1, a\rangle$ and then $xy \in D_G (1,a) = 1 + a$ $\equiv_G$.
  \item Assume that exist $p \in \dot M(G)$ with
  
   \begin{align*}
   a & \in c + cp \\
   b & \in p + ap \\
   d & \in p + cp.
  \end{align*}
  
   Then $a \in D(c,cp), b \in D(p,ap), d \in D(p,cp)$ and so 

  \begin{align*}
  \langle a, ap \rangle & \equiv \langle c, cp\rangle \\
  \langle b, ab \rangle & \equiv \langle ap, p \rangle \\
  \langle d, cd \rangle & \equiv \langle cp, p\rangle.
  \end{align*}

  Therefore $\langle a,b,ab \rangle \equiv \langle c, d, cd \rangle$. But since $G$ is special group, $\langle a,b,ab\rangle \equiv \langle 
d,c,cd\rangle$. Thus exists 
  $x,y,z \in G$ such that

  \begin{align*}
  \langle a, x \rangle & \equiv \langle d, y\rangle \\
  \langle b, ab \rangle & \equiv \langle x, z \rangle \\
  \langle c, cd \rangle & \equiv \langle y, z\rangle.
  \end{align*}
  
  It is possible to conclude that $a = b (ab) = xz$ and $d = yz$. So $x = az, y = dz$ and then
  
  \begin{align*}
  \langle a, az \rangle & \equiv \langle d, dz\rangle \\
  \langle b, ab \rangle & \equiv \langle az, z \rangle \\
  \langle c, cd \rangle & \equiv \langle dz, z\rangle.
  \end{align*}
  
  We conclude $a \in d + dz, b \in z + az, c \in z + dz$ as desired.
 \end{enumerate}
\end{proof}

\begin{defn}\label{defn:special.mf}
 A multifield $F$ satisfying the properties i-iv of proposition \ref{prop:missues} will be called 
a \textbf{special multifield}. Note that, if $G$ is a special group, then $M(G)$ is a special multifield.
\end{defn}

\begin{prop}
 Every real reduced multifield is a special multifield.
\end{prop}
\begin{proof}
 Let $F$ be a Real Reduced Multifield. If $a\in\dot F$, from $a^3=a$, cancelling we obtain $a^2=1$.
\begin{align*}
 (1+a)(1+a)\subseteq1+a+a+a^2=(1+1)+(a+a)=1+a.
\end{align*}

Let $x\in1+a$. We get
\begin{align*}
 x\in1+a&\Rightarrow1=x^2\in x+ax\Rightarrow ax\in1+a.
\end{align*}
Then $1=(ax)(ax)\in(1+a)(1+a)=1+a$.

Now, let $a,b,c, p \in \dot{F}$ such that
\begin{align} \label{34}
   a & \in c + cp \\
   b & \in p + ap \\
   d & \in p + cp.
  \end{align}
  
  Multiplying the last equation by $pcd$ we obtain $cp \in d + cd$. Then $a \in c + (d + cd)$. By the associativity property of $F$, exist $l' 
\in c + cd$ with $a \in d + l'$.
  If $l' = 0$, then $d = -1, a = d$ and $l = p$ satisfies $a \in d + ld, b \in l + al, c \in l + dl$. If $l' \neq 0$, then let $l = l'd$. Thus
  
  \begin{align} \label{35}
   a & \in d + dl \\
   c & \in l + dl.
  \end{align}
  
  To complete the proof, we need to show $b \in l + al$. Let $\sigma \in \mbox{sper}(F)$. If $\sigma(a) = -1$, then $\sigma(b) \in \sigma(l) + 
\sigma(al)$ (*).
  If $\sigma(a) = 1$, then we have two more cases. If $\sigma(l) = 1$, then $\ref{35}$ implies $\sigma(d) = 1 = \sigma(c)$ and so $\ref{34}$ 
gives $\sigma(b) = 1$. This proves
  $(*)$. If $\sigma(l) = -1$, then by similar case analysis we conclude $\sigma(b) = -1$ and so $(*)$ is verified. Then by the local-global 
principle $\ref{cor:4.4marshall}$
  $b \in l + al$.
\end{proof}

\begin{teo}\label{mf.to.special}
 If $F$ is a special multifield the $(\dot F,\equiv,-1)$ is a special group where $\langle 
a,b\rangle\equiv\langle c,d\rangle\Leftrightarrow ab=cd$ and $a\in c+d$.
\end{teo}
\begin{proof}
 By (i), we have $(\dot F,1)$ is a group of exponent 2. Now, we will check each axiom of 
definition \ref{defn:sg}:
 \begin{description}
  \item [SG0] - By (ii) $1\in 1 + ab$, so $ab\in 1+ab$ and $a\in b+a$. As $ab=ab$, then 
$\langle a,b\rangle\equiv\langle a,b\rangle$, i.e,  $\equiv$ is reflexive. If $\langle 
a,b\rangle\equiv\langle c,d\rangle$, then
$ab=cd$ and $a\in c+d$. Then $ab\in cb + db$, so by $ab = cd$, we have $cd \in ad + db$ and then $c \in a + 
b$. So $\langle c,d\rangle\equiv\langle a,b\rangle$ and 
$\equiv$ is symmetric. Finally, suppose that $\langle a,b\rangle\equiv\langle c,d\rangle$ and $\langle 
c,d\rangle\equiv\langle e,f\rangle$. First, $ab=cd$ and 
$cd=ef$ implies $ab=ef$. Second, in order to show that $a \in e + f$, note that $a\in c+d \Rightarrow ac \in 1 
+ cd = 1 + ef$ and $c\in e+f \Rightarrow ce \in 1 + ef$; then by (iii), we have $ae \in 1+ef$ and so $a \in e 
+ f$. Therefore 
$\langle a,b\rangle\equiv\langle e,f\rangle$.
  
  \item [SG1] - As $F$ is a multifield, $ab=ba$. By (ii), $1\in 1 + ab$, then $ab\in 1+ba$ and $b\in a+b$. 
Therefore $\langle a,b\rangle\equiv\langle b,a\rangle$. 

  \item [SG2] - Since $1 \in 1 - a$, we have $a \in 1 - 1$. Therefore $\langle 
a,-a\rangle\equiv\langle1,-1\rangle$.

  \item [SG3] - Follows by definition.
  
  \item [SG4] - $\langle a,b\rangle\equiv\langle c,d\rangle\Rightarrow ab=cd$ and $a\in c+d$.
  \begin{align}\label{eq:1}
  ab=cd\Rightarrow -abbc=-bccd\Rightarrow-ac=-bd 
  \end{align}
  \begin{align}\label{eq:2}
   a\in c+d\Rightarrow ad\in 1+cd = 1 + ab\Rightarrow d\in a+b\Rightarrow a\in-b+d
  \end{align}
  so by \ref{eq:1} and \ref{eq:2} follow that $\langle a,-c\rangle\equiv\langle -b,d\rangle$.
  
  \item [SG5] - $\langle a,b\rangle\equiv\langle c,d\rangle\Rightarrow ab=cd$ and $a\in c+d$ 
$\stackrel{I}{\Rightarrow}$ $(ga)(gb)=(gc)(gd)$ and $ga\in gc+gd$ $\Rightarrow$ $\langle 
ga,gb\rangle\equiv\langle gc,gd\rangle$.

  \item [SG9] - The proof is analogous to $\ref{prop:missues}$, item $iv$.
  
  \end{description}
  
\end{proof}

\begin{cor}\label{cor:equiv2}
 In the objects of $\mathcal{SMF}$, define $S(F)=\dot F$ as the special group as stated in theorem \ref{mf.to.special}. Now, let 
$\sigma:F\rightarrow K$ be a SMF-morphism and define $S(\sigma)=\sigma|_{\dot F}$. Then $S:\mathcal{SMF}\rightarrow\mathcal{SG}$ 
is a functor.
\end{cor}
\begin{proof}
  We have $S(\sigma)$ is a group homomorphism with $S(\sigma)(-1)=-1$. If $a,b\ne0$ and $c\in a+b$, $c\ne0$, then there 
exists $d\in \dot F$ such that $\langle a,b\rangle\equiv_{S(F)}\langle c,d\rangle$, and as $c\in a+b\rightarrow 
\sigma(c)\in\sigma(a)+\sigma(b)$, we have $\langle\sigma(a),\sigma(b)\rangle\equiv_{S(K)}\langle 
\sigma(c),\sigma(d)\rangle$. Therefore: 
\begin{align*}
 (c\in a+b\rightarrow\sigma(c)\in\sigma(a)+\sigma(b))\Rightarrow \\
 (c\in D_{S(F)}(a,b)\rightarrow\sigma(c)\in D_{S(K)}(\sigma(a),\sigma(b)))
\end{align*}
And $S(\sigma)$ is a SG-morphism. Applying the same argument, we proof that $S(\sigma\tau)=S(\sigma)S(\tau)$. Hence, $S$ is a 
morphism.
\end{proof}

\begin{teo}\label{teo:sgsmfequiv}
 $S:\mathcal{SMF}\rightarrow\mathcal{SG}$ gives an equivalence of categories between $\mathcal{SG}$ and $\mathcal{SMF}$.
\end{teo}
\begin{proof}
 By the corollaries \ref{cor:equiv1} and \ref{cor:equiv2}, we have functors 
$M:\mathcal{SG}\rightarrow\mathcal{SMF}$ and $S:\mathcal{SMF}\rightarrow\mathcal{SG}$. We will proof 
that $M\circ S\cong Id_{\mathcal{SMF}}$ and $S\circ M\cong Id_{\mathcal{SG}}$.
\begin{enumerate}[i -]
 \item $M\circ S\cong Id_{\mathcal{SMF}}$. Let $F$ be a special multifield. How $S(F)=\dot F$ and 
$M(S(F))=S(F)\cup\{0\}$, we have $M(S(F))=F$. Next, let $\sigma:F\rightarrow K$ be a SMF-morphism. We 
have $S(\sigma)=\sigma|_{\dot F}$ and $M(S(\sigma))$ is defined with the extension 
$S(\sigma)(0)=0$. Therefore $M(S(\sigma))=\sigma$ and $M\circ S\cong Id_{\mathcal{SMF}}$.
 
 \item $S\circ M\cong Id_{\mathcal{SG}}$. Let $G$ be a special group. Again, $M(G)=G\cup\{0\}$ and 
$S(M(G))=M(G)\setminus\{0\}$. Hence $S(M(G))=G$. Next, let $f:G\rightarrow H$ be a SG-morphism. How 
$M(f)$ is defined with the extension $f(0)=0$ and $S(M(f))=M(f)|_{M(G)\setminus\{0\}}$, we have 
$S(M(f))=f$ and $S\circ M\cong Id_{\mathcal{SG}}$, finalizing the proof.
\end{enumerate}
\end{proof}

$ $

We can summarize the functors obtained by the following diagram:
$$\xymatrix{ & \mathcal{AOS}^{op}\ar@<.9ex>[dr]_{\cong}\ar@<.9ex>[dl] & \\
\mathcal{RSG}\ar@<.9ex>[ur]_{\cong}\ar@<.9ex>[rr]_{\cong}\ar@<-.9ex>@{^{(}->}[d] & & 
\mathcal{MF}_{red}\ar@<.9ex>[ll]\ar@<.9ex>[lu]\ar@<-.9ex>@{^{(}->}[d] \\
\mathcal{SG}\ar@<.9ex>[rr]_{\cong} & & \mathcal\mathcal{SMF}\ar@<.9ex>[ll]}$$

\begin{teo}
 Let $M:\mathcal{SG}\rightarrow\mathcal{SMF}$ the functor defined in \ref{cor:equiv1}.
 \begin{enumerate}[i -]
  \item $M$ preserves products.
  \item $M$ preserves quotients.
  \item $M$ preserves directed limits. 
 \end{enumerate}
\end{teo}
\begin{proof}
 $ $
\begin{enumerate}[i -]
 \item Firstly, observe that $\mathcal{SMF}$ has products, because the categorical equivalence with 
$\mathcal{SG}$. However, this product is not the restriction of the product in $\mathcal{MF}$.

Now, let $\{G_i\}_{i\in I}$ be a family of special groups. The product $G=\prod^n_{i\in I}G_i$ is defined 
with the operation and special relation given pontwise, and $-1=(-1,-1,...)$, i.e,
$$\langle(a_i)_{i\in I},(b_i)_{i\in I}\rangle\equiv_G\langle(c_i)_{i\in I},(d_i)_{i\in I}\rangle
\Leftrightarrow\langle a_i,b_i\rangle\equiv_{G_i}\langle c_i,d_i\rangle,\,\forall\,i\in I.$$
This implies that $(a_i)_{i\in I}D_G((c_i)_{i\in I},(d_i)_{i\in I})$ if and only if
$a_i\in D_{G_i}(c_i,d_i)$ for all $i\in I$. This argument shows that
$$M\left(\prod^n_{i\in I}G_i\right)=\prod^n_{i\in I}M(G_i).$$

 \item More specifically, we want to show that if $G$ is a special group and $\Delta\subseteq G$ is a saturated 
subgroup
\footnote{We say that $\Delta$ is \textbf{saturated} if for all $a\in G$, $a\in\Delta\Rightarrow 
D_G(1,a)\subseteq\Delta$.} then $M(G/\Delta)\cong M(G)/\tilde{\Delta}$, when 
$\tilde{\Delta}=\{M(\delta):\delta\in\Delta\}$. The isometry relation on the 
quotient group $G/\Delta$ is:
$$\langle a/\Delta,b/\Delta\rangle\equiv^*_G\langle c/\Delta,d/\Delta\rangle\mbox{ if and only if }
\begin{cases}
 \exists\,a',b',c',d'\in G\mbox{ such that } \\
 aa',bb',cc',dd'\in\Delta\mbox{ and }\\
 \langle a',b'\rangle\equiv_G\langle c',d'\rangle.
\end{cases}$$
 This implies that $a/\Delta\in D_{G/\Delta}(c/\Delta,d/\Delta)$ if and only if there exist $r,s,t\in G$ such that
 $r\in D_G(s,t)$, with $ar,cs,dt\in\Delta$. Multiplying this by $arcsdt\in\Delta$, we have
 $a(csdt)\in D_G(c(ardt),d(arcs))$, and $csdt,ardt,arcs\in\Delta$. Applying the functor, we have $\overline 
a\in\overline c+\overline d$ in $M(G)/\tilde{\Delta}$, and the desired follow by this.
 
 \item Let $\mathcal G=(G_i,\{f_{ij}:i\le j\},I)$ be an inductive system of special groups. Let $G$ be the 
inductive limit of $\mathcal G$ and let $f_i:G_i\rightarrow G$ the correspondent SG-morphism associated to 
this construction. Then given $\langle a,b\rangle\equiv_G\langle c,d\rangle$ if and only if there exist $i\in I$ and 
$a_i,b_i,c_i,d_i\in G_i$ such that $\langle a_i,b_i\rangle\equiv_{G_i}\langle c_i,d_i\rangle$ and
$\langle f_i(a_i),f_i(b_i)\rangle=\langle a,b\rangle$, $\langle f_i(c_i),f_i(d_i)\rangle=\langle 
c,d\rangle$ (both over $G$). This is suffice to show that
$$M\left(\varinjlim\limits_{i\in I}G_i\right)=\varinjlim\limits_{i\in I}M(G_i).$$
\end{enumerate}
\end{proof}

\section{Real Reduced Multirings, Abstract Real Spectra and Real Semigroups}

  Since abstract ordering spaces and special groups generalizes almost entire classical and reduced theory of quadratic forms 
over fields, we could (naturally) ask the following
\begin{center}
 \textit{Is there some reasonable theory of quadratic forms over general coefficients in rings?}
\end{center}

There is an excellent book, \cite{knus1991quadratic}, that deal with quadratic forms in an style near to that was presented in 
Lam's classical books \cite{lam1983orderings} and \cite{lam2005introduction}, in the most general possible setting. And of course, 
some abstract theories appears trying to deal with this question. In 90's Marshall generalizes the AOS to rings, and 
called his new theory by ``Abstract Real Spectrum''. But the ring-theoretic case is much more difficult that the 
field one, the isometry is not well behaved and an algebraic counterpart of the ARS's appears just in years 2000, with the 
real semigroups (RS) of Dickmann and Petrovich \cite{dickmann2004real}. 

The RS appears in an atempt to creat a duality $\mathcal{RS}\simeq\mathcal{ARS}^{op}$ likewise 
$\mathcal{SG}\simeq\mathcal{AOS}^{op}$, goal that was successfully achieved. This theory is still in development: there is a 
preliminary book \cite{dickmann2012real}, which cover almost of the basics aspects of the theory. From the middle of the years 
2010's, Dickmann and Petrovich expanded the theory in \cite{dickmann2012spectral} and with the participation of F. Miraglia, 
more sophisticate constructions appears firstly in \cite{dickmann2017constructions}. Dickmann and Miraglia achieved another step 
in \cite{dickmann2015faithfully}: they exhibit a new description of the invertible elements of a real semigroup associated to a 
ring with many units.

In this section, we will construct (again) two equivalence of categories. The first one is given by a functor 
$\mathcal{MR}_{red}\rightarrow\mathcal{ARS}^{op}$ where $\mathcal{ARS}$ is the category of abstract real spectra as in 
\cite{marshall1996spaces} and $\mathcal{MR}_{red}$ is the category of real reduced multirings. This functor was indicated by 
Marshall in \cite{marshall2006real}, so we will establish precisely its definition and examine some properties. This encoding 
reveals how to extend the category of real semigroups to formally real semigroup (FRS), here introduced. The second equivalence 
extends the first one, that is alternatively given by a functor $\mathcal{MR}_{red}\rightarrow\mathcal{RS}$, from the category of 
real reduced multirings to the category of real semigroups. This second equivalence provides a 
functor $M:\mathcal{FRS}\hookrightarrow\mathcal{MR}$, here introduced, when $\mathcal{MR}$ is the category of multirings. The 
image of this functor is a subcategory of $\mathcal{MR}$, that we will call \textbf{special multirings}, and denote by 
$\mathcal{SMR}$.
 
Recall that, $\{-1,0,1\}$ has a natural ordering relation and for any set $X$, $\{-1,0,1\}^X$ denotes 
the set of all functions $a:X\rightarrow\{-1,0,1\}$. This is a monoid with the operation given by 
$(ab)(x)=a(x)b(x)$. 
 
\begin{defn}[Abstract Real Spectra]\label{defn:ars}
 An \textbf{abstract real spectra or space of signs}, is a pair $(X,G)$ satisfying:
 \begin{description}
  \item [AX1 -] $X$ is a non-empty set, $G$ is a submonoid of $\{-1,0,1\}^X$, $G$ contains the constants 
functions $-1,0,1$, and $G$ separates points in $X$. 
 \end{description}

 If $a,b\in G$, the \textbf{value set} $D(a,b)$ is defined to be the set of all $c\in G$ such that, for 
all $x\in X$, either $a(x)c(x)>0$ or $b(x)c(x)>0$ or $c(x)=0$. The \textbf{value set} $D^t(a,b)$ is 
defined to be the set of all $c\in G$ such that, for all $x\in X$, either $a(x)c(x)>0$ or $b(x)c(x)>0$ 
or $c(x)=0$ and $b(x)=-a(x)$. Note that $c\in D^t(a,b)\Rightarrow c\in D(a,b)$. Conversely, $c\in 
D(a,b)\Rightarrow c\in D^t(ac^2,bc^2)$.

\begin{description}
 \item [AX2 -] If $P$ is a submonoid of $G$ satisfying $P\cup-P=G$, $-1\notin P$, $a,b\in P\Rightarrow 
D(a,b)\subseteq P$ and $ab\in P\cap-P\Rightarrow a\in P\cap-P$ or $b\in P\cap-P$, then there exists 
$x\in X$ (necessarily unique) such that $P=\{a\in G:a(x)\le0\}$.

\item [AX3 (Strong Associativity) -] For all $a,b,c\in G$, if $p\in D^t(a,r)$ for some $q\in D^t(b,c)$ 
then $p\in D^t(r,c)$ for some $r\in D^t(a,b)$.
\end{description}
\end{defn}
 
\begin{defn}\label{defn:arsmorphism}
 A \textbf{morphism} of abstract real spectras $(X,G)\rightarrow(Y,H)$, or an \textbf{ARS-morphism}, is a mapping 
$\tau:X\rightarrow Y$ such that for 
each $a\in H$, the composite mapping is $a\circ\tau:X\rightarrow\{-1,0,1\}$ is an element of $G$ (so 
$\tau$ is surjective and induces a mapping $a\mapsto a\circ\tau$ from $H$ to $G$). $\tau$ is said to be 
an \textbf{isomorphism} if the mappings $X\rightarrow Y$ and $H\rightarrow G$ are bijective.
\end{defn}

\begin{teo}\label{teo:arstomrred}
 Let $(X,G)$ an abstract real spectra and define $a+b= D^t(a,b)$. Then 
$(G,+,\cdot,-,0,1)$ is a real reduced multiring.
\end{teo}
\begin{proof}
 Firstly, observe that $+$ is well-defined. Then, we will verify the conditions of 
definition \ref{defn:multiring}. Commutativity, associativity and neutral element ($a\in 
D^t(0,b)\Leftrightarrow a=b$) are immediate. In fact, the unique non-trivial part of the proof is
$$a\in D^t(b,c)\Rightarrow b\in D^t(a,-c)\mbox{ and }c\in D^t(-b,a).$$
We will prove that $b\in D^t(a,-c)$ and the case $c\in D^t(-b,a)$ analogous. Let $x\in X$ and $a\in 
D^t(b,c)$. Remember that $a\in D^t(b,c)$ means that $a(x)b(x)>0$ or $a(x)c(x)>0$ or $a(x)=0$ and 
$b(x)=c(x)$ happens for all $x\in X$.

If $a(x)b(x)>0$, then $b(x)a(x)>0$ and it is done. If $a(x)c(x)>0$, we have some cases:
\begin{itemize}
 \item $a(x)=c(x)=1$. We can suppose that $a(x)b(x)\le0$ and $b(x)\in\{0,1\}$. If $b(x)=0$ it is done. 
If $b(x)=1$, then $b(x)[-c(x)]>0$.
 
 \item  $a(x)=c(x)=1$. Again, we will suppose that $a(x)b(x)\le0$ and $b(x)\in\{0,1\}$. If $b(x)=0$ it 
is done. If $b(x)=1$, then $b(x)[-c(x)]>0$.

\item $a(x)=0$ and $b(x)=c(x)$. If $b(x)=c(x)=0$ then $b(x)=0$ and $a(x)=c(x)$. If $b(x)=c(x)\ne0$, then 
$b(x)c(x)>0$.
\end{itemize}
Hence $G$ is a multiring. For the real reduced part, we have immediately that $1\ne0$ and $a^3=a$ for 
all $a\in G$.
\begin{align*}
 c\in D^t(a,ab^2)\Leftrightarrow \forall x \in X c(x)a(x)=0\vee (c(x)=0\wedge a(x)=0)\Leftrightarrow c=a
\end{align*}
and 
\begin{align*}
 c&\in D^t(a^2,b^2)\Leftrightarrow \\
 \forall\,x&\in X ((c(x)=1 \ \wedge (a^2(x) = 1 \vee b^2(x) = 1) ) \vee(c(x)=0\wedge a(x)b(x)=0))
\end{align*}
This implies that $c$ is uniquely determined. Therefore, $G$ is a real reduced multiring.
\end{proof}

\begin{cor}\label{cor:equiv6}
 Let $(X,G)$ and $(Y,H)$ be two abstract real spectra and $\tau:Y\rightarrow X$ be an ARS-morphism. Define 
$M(X)$ as the real reduced multiring as in theorem \ref{teo:arstomrred} and $M(\tau)=f$ when 
$f:G\rightarrow H$ is the group homomorphism induced by $\tau$. Then  $M:\mathcal{ARS}^{op}\rightarrow\mathcal{MR}_{\mbox{red}}$ 
is a functor.
\end{cor}
\begin{proof}
  We have $c\in a+b\Rightarrow c\in D^t(a,b)\Rightarrow f(c)\in D^t(f(a),f(b))\Rightarrow f(c)\in f(a)+f(b)$ by an argument 
analogous to the corollary \ref{cor:equiv5}. Then $M(\tau)$ is a multiring morphism and this is suffice to prove that $M$ 
is a (contravariant) functor.
\end{proof}

\begin{teo}\label{teo:mrredtoars}
 Let $A$ be an real reduced multiring and consider the strong embedding 
 $$i:A\rightarrow Q_2^{\mbox{Sper}(A)}$$ 
 given by $i(a)=\hat a:\mbox{Sper}(A)\rightarrow Q_2$ when $\hat a(\sigma)=\sigma (a)$. Define $\hat A=i(A)$. Then $(\mbox{Sper}(A),\hat 
A)$ is an abstract real spectra. 
\end{teo}
\begin{proof}
 We will check each definition of \ref{defn:ars}:
 \begin{description}
  \item [AX1 -] Is consequence of the multiring structure on $\hat A$ (with the fact that $\hat A$ is contained in $Q_2^{\mbox{Sper}(A)}$).
  
  \item [AX2 -] Let $P$ be a submonoid of $\hat A$ such that $P\cup-P=\hat A$, $-1\notin P$, 
$a,b\in P\Rightarrow D(a,b)\subseteq P$ and $ab\in P\cap-P\Rightarrow a\in P\cap-P$ or $b\in P\cap-P$. 
Firstly, observe that
\begin{align}\label{eq:foratemer}
 D^t(a,b)=\{d:d\in a+b\}.
\end{align}
In fact, $d\in D^t(a,b)$ if and only if $\forall\sigma\in\mbox{Sper}(A)$, $\sigma(d)\sigma(a)>0$ or 
$\sigma(d)\sigma(b)>0$ or $\sigma(d)=0$, and $\sigma(a)=-\sigma(b)$ if and only if 
$\sigma(d)\in\sigma(a)+\sigma(b)$ for all $\sigma\in\mbox{Sper}(A)$. By the local-global principle 
for multirings \ref{prop:7.3marshall} we have this happens if and only if $d\in a+b$.
  
  \item [AX3 -] This is consequence of \ref{eq:foratemer} and associativity.
 \end{description}
\end{proof}

\begin{teo}\label{teo:arsmrredequiv}
 Define $M:\mathcal{ARS}^{op}\rightarrow\mathcal{MR}_{red}$ and $\mbox{Sper}:\mathcal{MR}_{red}\rightarrow\mathcal{ARS}^{op}$ as 
we already defined in corollary \ref{cor:equiv6} and theorem \ref{teo:mrredtoars}. These functors defines an equivalence of 
categories between $\mathcal{ARS}^{op}$ and $\mathcal{MR}_{\mbox{red}}$.
\end{teo}
\begin{proof}
  Follow  from $M\circ\mbox{Spec}\cong Id_{\mathcal{MR}_{red}}$ and $\mbox{Spec}\circ M\cong Id_{\mathcal{ARS}^{op}}$.
\end{proof}


\begin{defn}[Ternary Semigroup]\label{defn:tsemig}
 A \textbf{ternary semigroup} (abbreviated TS) is a structure $(S,\cdot,1,0,-1)$ with individual 
constants $1,,0,-1$ and a binary operation ``$\cdot$'' such that:
\begin{description}
 \item [TS1 -] $(S,\cdot,1)$ is a commutative semigroup with unity;
 \item [TS2 -] $x^3=x$ for all $x\in S$;
 \item [TS3 -] $-1\ne1$ and $(-1)(-1)=1$;
 \item [TS4 -] $x\cdot0=0$ for all $x\in S$;
 \item [TS5 -] For all $x\in S$, $x=-1\cdot x\Rightarrow x=0$.
\end{description}

We shall write $-x$ for $(-1)\cdot x$. The semigroup verifying conditions [TS1] and [TS2] (no extra 
constants) will be called \textbf{3-semigroups}. We denote $\mbox{Id}(S)=\{x\in S:x^2=x\}=S^2$ and $S^*=\{x\in S:x^2=1\}$.
\end{defn}

\begin{ex}
 $ $
 \begin{enumerate}[a -]
  \item The three-element structure $\bm 3=\{1,0,-1\}$ has an obvious ternary semigroup structure.
 \end{enumerate}
\end{ex}

Here, we will enrich the language $\{\cdot,1,0,-1\}$ with a ternary relation $D$. In agreement with 
\ref{defn:ars}, we shall write $a\in D(b,c)$ instead of $D(a,b,c)$. We also set:
$$a\in D^t(b,c)\Leftrightarrow a\in D(b,c)\wedge-b\in D(-a,c)\wedge -c\in D(b,-a).$$
The relations $D$ and $D^t$ are called \textbf{representation} and \textbf{transversal representation} 
respectively.

\begin{defn}[Real Semigroup]\label{defn:realsemigroup}
 A \textbf{real semigroup} is a ternary semigroup together with a ternary relation $D$ 
satisfying:
\begin{description}
 \item [RS0 -] $c\in D(a,b)$ if and only if $c\in D(b,a)$.
 \item [RS1 -] $a\in D(a,b)$.
 \item [RS2 -] $a\in D(b,c)$ implies $ad\in D(bd,cd)$.
 \item [RS3 (Strong Associativity) -] If $a\in D^t(b,c)$ and $c\in D^t(d,e)$, then there exists $x\in 
D^t(b,d)$ such that $a\in D^t(x,e)$.
 \item [RS4 -] $e\in D(c^2a,d^2b)$ implies $e\in D(a,b)$.
 \item [RS5 -] If $ad=bd$, $ae=be$ and $c\in D(d,e)$, then $ac=bc$.
 \item [RS6 -] $c\in D(a,b)$ implies $c\in D^t(c^2a,c^2b)$.
 \item [RS7 (Reduction) -] $D^t(a,-b)\cap D^t(b,-a)\ne\emptyset$ implies $a=b$.
 \item [RS8 -] $a\in D(b,c)$ implies $a^2\in D(b^2,c^2)$.
\end{description}
\end{defn}

The theory of real semigroups can be alternatively axiomatized by the transversal relation $D^t$. In this case, we define
$$c\in D(a,b)\Leftrightarrow c\in D^t(c^2a,c^2b).$$

The definition of morphism is quite standard: $f:(G,\cdot,1,0-1)\rightarrow(H,\cdot,1,0-1)$ is an 
\textbf{RS-morphism} if $f:G\rightarrow H$ is a morphism of semigroups, (i.e, $f(ab)=f(a)f(b)$, $f(1)=1$ 
and $f(0)=0$); $f(-1)=-1$ and $a\in D(b,c)\Rightarrow f(a)\in D(f(b),f(c))$ (hence $a\in 
D^t(b,c)\Rightarrow f(a)\in D^t(f(b),f(c))$). The category of real semigroups and their morphisms will be 
denoted by $\mathcal{RS}$.

A \textbf{formally real semigroup} is a ternary semigroup together with a ternary relation $D$ 
satisfying [RS0]-[RS6], [RS8] and:
\begin{description}
 \item [RS7a (Zero) -] $D^t(0,a)=\{a\}$.
 \item [RS7b (Semi-reality) -] For all $n\ge1$, $a_1,...,a_n\in G$, $-1\notin D^t(a_1^2,...,a_n^2)$, with the conventions 
$D^t(a)=\{a\}$ and
 $$D^t(a_1,...,a_n):=\bigcup\limits_{c\in D^t(a_2,...,a_n)}D^t(a,c).$$
\end{description}

The definition of morphisms of a formally real semigroup is analogous. The category of formally 
real 
semigroups and their morphisms will be denoted by $\mathcal{FRS}$.

\begin{prop}\label{prop:rsprop}
 The properties below holds in any formally real semigroup $G$, for all $a,b,c,d\in G$:
 \begin{enumerate}
  \item $a\in D(b,c)\Leftrightarrow a\in D^t(a^2b,a^2c)$.
  \item $a\in D^t(b,c)\Rightarrow -b\in D^t(-a,c)$.
  \item $0\in D(a,b)$.
  \item $a\in D^t(b,c)\Rightarrow ad\in D^t(bd,cd)$.
  \item $d\in D(ca,cb)\Rightarrow d=c^2d$. In particular, $D(0,a)\subseteq\{a^2x:x\in G\}$.
  \item $a^2\in D(1,b)$.
  \item $a\in D(0,0)\Leftrightarrow a=0$.
  \item $1\in D^t(1,a)$.
  \item $D^t(1,-1)=G$.
  \item $ab\in D(1,-a^2)$.
  \item $D^t(a,b)\ne\emptyset$.
  \item (Weak Associativity) $a\in D(b,c)\wedge c\in D(d,e)\Rightarrow\exists\,x[x\in D(b,d)\wedge a\in 
D(x,e)]$.
 \end{enumerate}
 Moreover, if $G$ is a real semigroup, then: 
  \begin{enumerate}
  \setcounter{enumi}{12}
  \item For all $a,b\in G$, $D^t(a^2,b^2)=\{x^2\}$ for some $x\in G$.
  \item $0\in D^t(a,b)\Leftrightarrow a=-b$ is equivalent to RS7a.
  \item $a\in D(0,1)\cup D(1,1)\Rightarrow a=a^2$.
  \item $a\in D^t(b,b)\Leftrightarrow a=b$.
 \end{enumerate}
\end{prop}

\begin{prop}\label{cor:3dickman}
 The ternary semigroup $\bm3=\{1,0,-1\}$ has a unique structure of real semigroup, with representation given by:
$$\begin{cases}
   D_{\bm3}(0,0)=\{0\}; \\
   D_{\bm3}(0,1)=D_{\bm3}(1,0)=D_{\bm3}(1,1)=\{0,1\}; \\
   D_{\bm3}(0,-1)=D_{\bm3}(-1,0)=D_{\bm3}(-1,-1)=\{0,-1\}; \\
   D_{\bm3}(1,-1)=D_{\bm3}(-1,1)=\bm3
  \end{cases}$$
  and transversal representation given by:
$$\begin{cases}
   D^t_{\bm3}(0,0)=\{0\}; \\
   D^t_{\bm3}(0,1)=D^t_{\bm3}(1,0)=D^t_{\bm3}(1,1)=\{1\}; \\
   D^t_{\bm3}(0,-1)=D^t_{\bm3}(-1,0)=D^t_{\bm3}(-1,-1)=\{-1\}; \\
   D^t_{\bm3}(1,-1)=D^t_{\bm3}(-1,1)=\bm3\end{cases}$$  
\end{prop}
\begin{proof}
See corollary 2.4 in \cite{dickmann2004real}. 
\end{proof}

\begin{teo}[Separation Theorem from \cite{dickmann2004real}]\label{teo:septeo}
 Let $G$ be a RS, and $a,b,c\in G$ and $X_G=Hom(G,\bm3)$. Then:
 \begin{enumerate}[i -]
  \item $a\in D_G(b,c)$ if and only if for all $h\in X_G$, $h(a)\in D_{\bm3}(h(b),h(c))$.
  \item $a\in D^t_G(b,c)$ if and only if for all $h\in X_G$, $h(a)\in D^t_{\bm3}(h(b),h(c))$.
  \item If $a\ne b$, there is $h\in X_G$ such that $h(a)\ne h(b)$.
 \end{enumerate}
\end{teo}
\begin{proof}
 See theorem 4.4 in \cite{dickmann2004real}.
\end{proof}

\begin{cor}
 Every real semigroup is a formally real semigroup.
\end{cor}
\begin{proof}
 If $-1\in D^t(a^2,b^2)$, then for all $h\in X_G$, $-1=h(-1)\in D_{\bm3}(h(a^2),h(b^2))$. Since $h(a^2),h(b^2)\in\{0,1\}$ for all 
$h\in X_G$, $D_{\bm3}(h(a^2),h(b^2))\subseteq\{0,1\}$ for all $h\in X_G$, a contradiction. If $-1=x^2$, the same argument holds 
with $-1\in D^t(x^2,0)$.

By \ref{prop:rsprop}(13) the above argument is enough to show that $-1\notin D^t(a_1^2,...,a_n^2)$ for all $a_1,...,a_n\in G$, 
and all $n\ge1$.
\end{proof}

\begin{teo}\label{teo:rstomrred}
 Let $(G,\cdot,1,0,-1,D)$ be a real semigroup and define $+:G\times 
G\rightarrow\mathcal{P}(G)\setminus\{\emptyset\}$, $ a + b = D^t(a,b) $ and $-:G\rightarrow 
G$ by $-(g)=-1\cdot g$. Then $(G,+,\cdot,-,0,1)$ is a real reduced multiring.
\end{teo}
\begin{proof}
 Firstly, observe that by \ref{prop:rsprop}(xv) the sum is well-defined, i.e, $D^t(a,b)\ne\emptyset$ for 
all $a,b\in G$. 

Now, we will check that $G$ is a multiring: of course, by RS0 we have $a+b=b+a$ (i.e, 
$D^t(a,b)=D^t(b,a)$) and
\begin{align*}
 \begin{cases}
 d\in D^t(a,b)\Leftrightarrow d\in D(a,b)\wedge-a\in D(-d,b)\wedge -b\in D(a,-d) \\
 a\in D^t(d,-b)\Leftrightarrow a\in D(d,-b)\wedge -d\in D(-a,-b)\wedge b\in D(d,-a)\\
 b\in D^t(-a,d)\Leftrightarrow b\in D(-a,d)\wedge a\in D(-b,d)\wedge -d\in D(-a,-b)
 \end{cases}
\end{align*}

So $d\in D^t(a,b)\Rightarrow a\in D^t(d,-b)\wedge b\in D^t(-a,d)$, or in other words, $d\in 
a+b\Rightarrow a\in d+(-b)\wedge b\in(-a)+d$. If $x=y$, by RS1 $x\in 0+y$. Conversely, let $x\in0+y$. We 
just proved that $0\in x-y$ and $0\in y-x$ then by RS7, $x=y$. How RS3 states the associativity (like 
\ref{defn:multigroupII}) we have $G$ is a commutative multigroup.

Because the commutative semigroup structure of $(G,\cdot,-1,0,1)$, we have $(G,\cdot,1)$ is a 
commutative monoid and $a\cdot0=0$ for all $a\in G$. The distributive law is just 
\ref{prop:rsprop}(iii), we have $G$ is a multiring.

Finally, we prove that $G$ is real reduced. We already have $-1\ne0$ and $a^3=a$. We have too, that 
$1\in D^t(1,b^2)$ by \ref{prop:rsprop}(ix) then by 
\ref{prop:rsprop}(iii) $a\in D^t(a,ab^2)$. Now, how $t^3=t$ we have
\begin{align}\label{eq:treta}
 t\in D^t(v^2x,w^2y)\Rightarrow \nonumber \\
 t\in D(v^2x,w^2y)\wedge -v^2x\in D(-t^3,w^2y)\wedge-w^2y\in 
D(v^2x,-t^3) \nonumber \\
\stackrel{\mbox{RS4}}{\Rightarrow}t\in D(x,y)\wedge -v^2x\in D(-t,y)\wedge-w^2y\in D(x,-t)
\end{align}

Hence, how by RS1 $-a\in D(-a,-x)$ for all $a,x\in G$, follow
\begin{align*}
 x\in D^t(a,ab^2)\Rightarrow x\in D^t(a^2\cdot a,(ab)^2\cdot 
a)\stackrel{\ref{eq:treta}}{\Rightarrow} \\
x\in D(a,a)\wedge-a\in D(-x,a)\wedge-ab^2\in D(a,-x)\Rightarrow \\
[x\in D(a,a)\wedge-a\in D(-x,a)\wedge-a\in D(a,-x)]\wedge-ab^2\in D(a,-x)\Rightarrow \\
x\in D^t(a,a)\wedge-ab^2\in D(a,-x)\stackrel{\ref{prop:rsprop}(vii+x)}{\Rightarrow} x=a 
\end{align*}

Then $a+ab^2=\{a\}$. For the last property, we have by theorem \ref{teo:septeo}(ii), we have $d\in 
D^t(b^2,c^2)\Leftrightarrow h(d)\in D^t_{\bm3}(h(b^2),h(c^2))$ for every $h\in X_G$. How by proposition 
\ref{cor:3dickman} $D^t(t^2,s^2)$ is unitary for every $s,t\in\bm3$, we have $D^t(b^2,c^2)$ is 
unitary for every $b,c\in G$.

Hence, by definition \ref{defn:mrrealreduced} $G$ is a real reduced multiring.
\end{proof}

\begin{cor}\label{cor:equiv5}
 There is a functor $M:\mathcal{RS}\rightarrow\mathcal{MR}_{\mbox{red}}$.
\end{cor}
\begin{proof}
 Let $R,S\in\mathcal{RS}$ and $f:R\rightarrow S$ a RS-morphism. Define $M(R)$ how the real reduced 
multiring as in theorem \ref{teo:rstomrred} and $M(f)=f$. Of course, $M(f)$ is a multiring morphism, 
because $c\in a+b\Rightarrow c\in D^t(a,b)\Rightarrow f(c)\in D^t(f(a),f(b))\Rightarrow f(c)\in 
f(a)+f(b)$. This is suffice to prove that $M$ is a functor.
\end{proof}

In order to associate a real semigroup to each real reduced multiring, we are going to set down some facts about multirings:

\begin{prop}
   Let $A$ be a real reduced multiring. Then we have the following:
\begin{enumerate}[i -]
  \item $x\in a x^2 + b x^2$ if and only if $x \in a A^2 + b A^2$;
  \item $x\in a+b$ if and only if $x\in a x^2 + b x^2$, $-a \in b a^2 - x a^2$ and $-b\in a b^2 - x b^2$;
  \item If $ax = bx$, $ay = by$ and $z \in x z^2 + y z^2$, then $az = bz$;
  \item If $x \in a x^2 + b x^2$, then $x^2 \in a^2 x^2 + b^2 x^2$.
\end{enumerate}
\end{prop}
\begin{proof}
Since $A$ is a real reduced multiring, we have by the local-global principle for multirings \ref{prop:7.3marshall} that $a \in b + 
c$ if and only if $\sigma(a) \in \sigma (b) + \sigma (c)$ for all $\sigma \in \mbox{Sper}(A)$. So to prove these items we just 
need to do it in $Q_2$ which is trivial (it is just an amount of cases).
\end{proof}

\begin{teo}\label{teo:mrredtors}
 Let $A$ be a real reduced multiring. Then $(A,\cdot,1,0,-1,D)$ is a real semigroup, where $d\in D(a,b)\Leftrightarrow d\in 
d^2a+d^2b$.
\end{teo}
\begin{proof}
 Firstly, note that by the preceding proposition, $x \in D(a,b) \Leftrightarrow x \in a A^2 + b A^2$ and $D^t (a,b) = a + b$.
 
Now, we will check each axiom of definition 
\ref{defn:realsemigroup}:
 \begin{description}
  \item [RS0 -] Is just commutativity of sum.
  
  \item [RS1 -] It follows by item $i$ of the preceding proposition.
  
  \item [RS2 -] $a\in D(b,c)\Leftrightarrow a\in a^2b+a^2c \stackrel{d^3=d}{\Rightarrow} 
  ad\in (ad)^2bd+(ad)^2cd\Rightarrow ad\in D(bd,cd)$.
  
  \item [RS3 -] It is just associativity of the sum.
  
  \item [RS4 -] It follows by item $i$ of the preceding proposition.
  
  \item [RS5 -] It follows by item $iii$ of the preceding proposition.
  
  \item [RS6 -] It follows by the characterization of $D^t$.
  
  \item [RS7 -] Since in a real reduced multiring we have $a + a = a$, if exist $c \in a - b$ with $-c \in a-b$, then $0 \in c - c 
\in a - b + a -b = a - b$ and then $a = b$.
  
  \item [RS8 -] It follows by item $iv$ of the preceding proposition.
 \end{description}
\end{proof}

\begin{cor}\label{teo:rsmrredequiv}
 Define the functor $S:\mathcal{MR}_{\mbox{red}}\rightarrow\mathcal{RS}$ as in corollary \ref{cor:equiv5}. Then $S$ is an 
equivalence of categories between $\mathcal{RS}$ and $\mathcal{MR}_{\mbox{red}}$.
\end{cor}
\begin{proof}
 The proof of $S\circ M\cong Id_{\mathcal{RS}}$ and $M\circ S\cong Id_{\mathcal{MR}_{\mbox{red}}}$ is mutatis mutandis of theorem 
\ref{teo:sgsmfequiv}.
\end{proof}

Of course, we can adapt the proof of theorem \ref{teo:rstomrred} to obtain a functor 
$M:\mathcal{FRS}\hookrightarrow\mathcal{MR}$. The image of this functor is a subcategory of $\mathcal{MR}$, 
that we will call \textbf{special multirings}, and denoted by $\mathcal{SMR}$. Again, we can summarize the 
functors obtained by the following diagram:
$$\xymatrix{ & \mathcal{ARS}^{op}\ar@<.9ex>[dr]_{\cong}\ar@<.9ex>[dl] & \\
\mathcal{RS}\ar@<.9ex>[ur]_{\cong}\ar@<.9ex>[rr]_{\cong}\ar@<-.9ex>@{^{(}->}[d] & & 
\mathcal{MR}_{red}\ar@<.9ex>[ll]\ar@<.9ex>[lu]\ar@<-.9ex>@{^{(}->}[d] \\
\mathcal{FRS}\ar@<.9ex>[rr]_{\cong} & & \mathcal{SMR}\ar@<.9ex>[ll]}$$

\begin{cor}\label{propertiessgtom}
 Let $M:\mathcal{RS}\rightarrow\mathcal{MR}_{red}$ the functor defined in \ref{teo:rstomrred}. Then $M$ 
preserves products and directed limits.
\end{cor}
\begin{proof}
 Follow directly by the definition of product and directed limits in $\mathcal{RS}$.
\end{proof}

$ $

Finally, we provide a diagram for a better visualization of the functors obtained:

 \SelectTips{eu}{12}
 $$\xymatrix@!=2pc{
 \mathcal{RSG}\ar@{^{(}->}[dd]\ar@<.9ex>[rr]\ar[dr] &&
 \mathcal{AOS}^{op}\ar@{^{(}->}'[d][dd]\ar@<.9ex>[ll]\ar@<.9ex>[rr] &&
 \mathcal{MF}_{red}\ar@{^{(}->}[dd]\ar@<.9ex>[ll]\ar[dl]\\
 & \mathcal{SG}\ar@{^{(}->}[dd]\ar@<.9ex>[rr]
 && \mathcal{SMF}\ar@{^{(}->}[dd]\ar@<.9ex>[ll]
 \\
 \mathcal{RS}\ar@<.9ex>[rr]\ar[dr] &&
 \mathcal{ARS}^{op}\ar@<.9ex>[rr]\ar@<.9ex>[ll] &&
 \mathcal{MR}_{red}\ar@<.9ex>[ll]\ar[dl]\\
 & \mathcal{FRS}\ar@<.9ex>[rr]
 && \mathcal{SMR}\ar@<.9ex>[ll]
 }$$
 
\section{Conclusion and Future Works}

We have provided a functorial picture of categories of abstract theories of quadratic forms, connecting them with multifields and 
multirings. This brings new perspectives and methods to the abstract theories of quadratic forms in forthcoming papers, that will 
be briefly described in the sequel.

The generalization of the theory quadratic forms to general coefficients in rings is a hard step. The book 
\cite{knus1991quadratic} cover some basic aspects in the most general setting possible, and we have Marshall's theory of abstract 
real spectra (\cite{marshall1996spaces}) and its algebraic counterpart, the real semigroups of Dickmann and Petrovich 
(\cite{dickmann2004real}) given a nice approach for the \textbf{reduced} theory of quadratic forms on rings, but, most of the 
relevant aspects of quadratic forms, like Witt rings, Pfister forms and etc, are uncovered. In the 
forthcoming \cite{roberto2019nonreduced}, we propose the fundamentals for a \textbf{non reduced and first-order} abstract 
quadratic forms theory in general coefficients on rings, with the intuition and machinery of multirings and multifields, inspired 
by the functorial picture described here.

In the forthcoming \cite{ribeiro2019vonneuman}, we describe the reduction functor, presented in section 3, from the category of 
von Neumann regular hyperrings to the category of real reduced multirings as a \textbf{definable} functor. This is achieved 
through sheaf theoretic (scheme like) methods in the hyperring setting (again, inspired by our functorial picture). Applications 
of this will appear in others contexts, for example, more general abstract Witt rings in an attempt to obtain a more concrete 
theory of quadratic forms with general coefficients over rings, complementing \cite{roberto2019nonreduced}, which is not 
dependent on the signatures and orderings, that are not first-order definable.

{\bf Acknowledgements:} We want to express our profound admiration and sincere  gratitude to professors F. Miraglia and M. Dickmann. We also want 
 to thank the anonymous referee for her/his careful reading and valuable suggestions.

\bibliographystyle{plain}
\bibliography{multirings_quadratic_forms}
\end{document}